\documentclass[10pt]{article}

\usepackage{amsfonts}
\usepackage{amsmath,amssymb,amsthm}
\usepackage{url}
\usepackage{amsmath,amssymb,amsthm}\usepackage{amscd,color}
\newtheorem*{Tha}{Theorem A}
\newtheorem*{Thb}{Theorem B}
\newtheorem{Th}{Theorem}

\newtheorem{Lemma}{Lemma}

\newtheorem{Prob}{Problem}
\newtheorem*{LemmaA}{Lemma A}

\theoremstyle{remark}

\newtheorem{Rem}{Remark}
\newtheorem{rem}{Remark}
\newtheorem{Def}{Definition}
\newtheorem*{rema}{Remark A}
\newtheorem*{remb}{Remark B}
\newtheorem*{remc}{Remark C}
\newtheorem*{remd}{Remark D}
\newtheorem*{reme}{Remark E}

 \newcommand{\cotan}{\textrm{cotan}}

\newcommand{\p}{\partial}
\newcommand{\trace}{\mbox{\rm trace}}

\newcommand{\eup}{\mathfrak{p}}

\newcommand{\weg}[1]{}
\newcommand{\const}{\mbox{\rm const}}

 \def\registered{{\ooalign{\hfil\raise .00ex\hbox{\scriptsize R}\hfil\crcr\mathhexbox20D}}}

\begin{document}

\date{} 
\title{Two-dimensional  metrics admitting  precisely one  projective vector field\footnote{this article has an appendix of Alexei V. Bolsinov, 
Vladimir S. Matveev, and  
Giuseppe Pucacco}}
\author{
Vladimir S. Matveev\thanks{ Institute of Mathematics, FSU Jena, 07737 Jena Germany,  vladimir.matveev@uni-jena.de} \thanks{ partially supported by  DFG (SPP 1154  and GK   1523)}}

\maketitle\begin{abstract}
We give a complete list of  $2$-dimensional metrics
that admit an essential projective vector field.
 This solves a problem explicitly posed by Sophus Lie in 1882.
\end{abstract}
\section{Introduction}
\subsection{Main definitions  and results}\label{introduction}

Let $g$ be a smooth  Riemannian or pseudo-Riemannian metric
on a 2-dimensional disc $D^2$. \begin{Def} 
A vector field $v $ is called  \emph{projective}, if its flow takes (unparameterized) geodesics to geodesics.  \end{Def}  As Lie showed \cite{Lie},
the set of vector fields projective with respect to a given $g$
forms a Lie algebra (for our paper it is sufficient that it is  a  vector space). We will denote this Lie algebra by $\eup(g)$.

The following two problems  were  posed by Sophus Lie%
\footnote{
German original from \cite{Lie}, Abschn. I, Nr. 4, \\
{\bf Problem I:}
\emph{Es wird verlangt, die Form des Bogenelementes einer jeden Fl\"ache zu bestimmen,
deren geod\"atische Kurven eine  infinitesimale Transformation gestatten}.
\\ 
{\bf Problem II:}
\emph{Man soll die Form des Bogenelementes einer jeden Fl\"ache bestimmen,
deren geod\"atische Kurven mehrere infinitesimale Transformationen gestatten}.}
in 1882:

\begin{Prob}[Lie]
Find all metrics $g$ such that $\dim\left(\eup(g)\right)= 1$.
\end{Prob}

\begin{Prob}[Lie]
Find all metrics $g$ such that $\dim\left(\eup(g)\right)\ge 2$.
\end{Prob}

The second problem of Lie was completely solved in \cite{bryant}. The present paper 
 gives a solution  of the first problem of Lie. The reader should consult \cite{bryant,duna} for the history of the question, for the connection with the results of Aminova \cite{Aminova0,Aminova},  and for the description of the circle of ideas, though we recall some of them    in \S \ref{schema}.

The biggest family of metrics admitting projective vector fields consists of    metrics   admitting infinitesimal homotheties.    Recall that a vector field $v$ is  an 
\emph{infinitesimal homothety} for a metric $g$ if $L_vg = \lambda g$ for a certain constant $\lambda \in \mathbb{R}$, where $L_v$  denotes the Lie derivative. In this definition, we allow $\lambda  = 0$,   so that Killing vector fields are also infinitesimal homotheties.

This ``biggest" family  of metrics is  very well understood: it is well known and  it was explicitly mentioned by Lie in the paper  \cite{Lie},        that    in the coordinates  $(x,y)$ 
such that  $v =\frac{\partial }{\partial x}$ such a metric $g$ is given by the matrix 
 $e^{\lambda x}\begin{pmatrix} E(y) & F(y) \\ F(y) & G(y)   \end{pmatrix} 
 $, where $E,F,G$  are functions of $y$ only.

 Thus, the first Lie  Problem  as Sophus Lie himself understood it  is to find all 
 $g$ admitting no infinitesimal homotheties, but   such that  
 $\dim\left(\mathfrak{p}(g)\right)=1$.   
 From the paper \cite{Lie}  it is clear that Lie considered this problem only locally,  in a small neighborhood of a  generic point. 

The next three  theorems  solve   the Problem  1 (of Lie) above.

\begin{Def} 
Two metrics $g$ and $\bar g$ on $D^2$ are called  \emph{projectively equivalent}, if they have the same geodesics considered as unparameterized curves. 
\end{Def}

\begin{Th} \label{main} Assume the metric $\check g$ on $D^2$ admits a projective vector field $v$. Assume in addition that  for any open $U\subset D^2$ the restriction of $\check g$ to $U$  admits no infinitesimal homothety.

   Then, in a neighborhood of almost every point there exists a coordinate  
   system  $(x,y)$  (in certain cases  we  consider the corresponding  complex   coordinates $(z= x+ i\cdot y, \bar z= x-i\cdot y) $),     
   such that   in this neighborhood  the  vector field  $v$ and a certain metric 
    $g$ projectively equivalent to $\check g$  are given by the formulas below.    

    \begin{enumerate} 
    \item { \rm \bf (Liouville Case)}  $ds^2_{g} = (X(x)- Y(y))(X_1(x)dx^2+   Y_1(y)dy^2), v=  \frac{\partial }{\partial x} + \frac{\partial }{\partial y}$, where 

       \begin{enumerate} 
    \item  $X(x)= \frac{1}{x}$,  \label{1a}
    $Y(y)= \frac{1}{y}$, 
    $X_1(x)= c\cdot \frac{e^{-3x}}{x}$,
     $Y_1(y)=  \frac{e^{-3y}}{y}$. 
    \item \label{1b} $X(x)= \tan(x)$, $Y(y)= \tan(y)$, $X_1(x)= c\cdot \frac{e^{-3\lambda x}}{\cos(x)}$, $Y_1(y)= \frac{e^{-3\lambda y}}{\cos(y)}$. 
    \item \label{1c} $X(x)= c\cdot e^{\nu x}$, $Y(y)=  e^{\nu y}$, $X_1(x)=e^{2x}$, $Y_1(y)= \varepsilon  e^{2y}$.

    \end{enumerate} 
   \item  {\rm \bf (Complex-Liouville Case)}  $ds^2_{g} = (h(z)- \overline{h( z}))(h_1(z)dz^2- \overline{h_1( z})d\bar z^2),$ \  \  $v=   \frac{\partial }{\partial x} $ $(= \frac{\partial }{\partial  z} + \frac{\partial }{\partial \bar z})$,   where  \begin{enumerate} 
    \item \label{2a} $h(z)= \frac{1}{z}$, \  $h_1(z)= C\cdot \frac{e^{-3 z}}{z}$. 
    \item \label{2b} $h(z)= \tan(z)$,  $h_1(z)= C\cdot \frac{e^{-3\lambda z}}{\cos(z)}$. 
    \item  \label{2c} $h(z)= C \cdot e^{\nu z}$, $h_1(z)=e^{2z}$.

    \end{enumerate} 
   \item   {\rm \bf (Jordan-block  Case)  }   $ds^2_{g} = (Y(y) + x) dx dy $, $v= v_1(x,y)  \frac{\partial }{\partial x} + v_2(y) \frac{\partial }{\partial y}$, where  \begin{enumerate} 
    \item \label{3a}  $Y=  e^{\frac{3}{2y}}\cdot \frac{\sqrt{|y|}}{y-3}+  \int_{y_0}^y e^{\frac{3}{2\xi}}\cdot \frac{\sqrt{|\xi|}}{(\xi-3)^2}d\xi$, \\   $ v_1= \frac{y-3}{2} \left( x+ \int_{y_0}^y e^{\frac{3}{2\xi}}\cdot \frac{\sqrt{|\xi|}}{(\xi-3)^2}d\xi\right), $ 
    $v_2= y^2$.  

    \item \label{3b} $Y=  
    e^{-\frac{3}{2}\lambda \arctan(y)}\cdot \frac{\sqrt[4]{y^2+1}}{y-3\lambda}+ \int_{y_0}^y e^{-\frac{3}{2}\lambda \arctan(\xi)}\cdot \frac{\sqrt[4]{\xi^2+1}}{(\xi-3\lambda)^2}d\xi$, \\  
     $ v_1= \frac{y-3\lambda }{2} \left( x+  \int_{y_0}^y e^{-\frac{3}{2}\lambda \arctan(\xi)}\cdot \frac{\sqrt[4]{\xi^2+1}}{(\xi-3\lambda)^2}d\xi\right), $ $v_2= y^2+1$.

    \item \label{3c} $Y(y)= y^{\frac{1}{\eta}}$, $v_1(x, y)=  {x} $, $v_2={\eta} y$,    

      \item \label{3d} $Y(y)= y^{2}$, $v_1(x, y)= 2  x $, $v_2= y$,

    \end{enumerate}
    where $c \in \mathbb{R}\setminus \{0\}$, $y_0\in \mathbb{R}$, 
    $\lambda\in \mathbb{R}$, 
    $\nu, \eta \in  (0, 4]$, $\nu \ne 1$, $\eta\not\in \{\tfrac{1}{2}, 1\}$,  \   $C\in \mathbb{C}$, $|C|=1$, $\varepsilon \in \{-1, 1\}$  are constants, and $\overline{h}$ and $\overline{h_1}$ denote the complex-conjugate to $h$ and $h_1$. 

 Moreover,  in the case \ref{1b}, if $\lambda =0$, then $c\ne \pm 1$.  In the case \ref{2b}, if $\lambda =0$, then $C\ne \pm 1$. In  the case \ref{1c}, if $\nu = 2$, then $c\ne - \varepsilon$. In the case \ref{2c},  if $\nu =2$, then $C\ne\pm1 $.     

    \end{enumerate}

\end{Th}

\begin{Rem} We do not claim in Theorem \ref{main} that all   metrics   projectively equivalent to $g$ admit no infinitesimal homotheties.   In view of  Theorem \ref{main3}, it is easy to understand  whether a metric $\bar g$
 projectively equivalent to $g$ from Theorem \ref{main} admits an infinitesimal  homothety: indeed, by Theorem \ref{main3},  the metrics from Theorem \ref{main} have unique (up to multiplication by a constant)  projective vector  field. Thus it is sufficient to  check whether $v$ is an infinitesimal  homothety.

 Moreover, from  the proof of Theorems \ref{main}, \ref{main2}  it will be clear that for every metric $ g$ from Theorem \ref{main}  the set of metrics projectively equivalent to $g$ and admitting an infinitesimal homothety is very small (has dimension at most  1 in the two- or three-dimensional space of metrics projectively equivalent to $g$.) 

 But certain metrics projectively equivalent to $g$ may admit infinitesimal homotheties. For example, in the case  \ref{1c},   the vector field $v $  is already an infinitesimal 
 homothety for  $g$.
\end{Rem}

Clearly, projective equivalence is a symmetric, reflexive and transitive relation on the space of all metrics  on $U\subseteq D^2$, i.e., it is an equivalence relation. The equivalence class of a metric $g$ with respect to projective equivalence will be  called the \emph{projective class}  of a metric 
(we give an equivalent analytic  definition in \S\ref{2.1}), it contains all metrics projectively equivalent to $g$.  
 Clearly, if $v$ is a projective vector field for  a metric from a projective class, it is a projective vector field for every metric from the projective class.  Theorem \ref{main} describes   (in a neighborhood of almost every point) all  projective classes admitting essential  projective vector fields. The next theorem describes all metrics of such projective classes. 

 For two   metrics (three metrics, respectively)  
  $g$ and $\bar g$ on $U\subseteq D^2$  ($g$, $\bar g$,  and  $\tilde g$, respectively) 
   and for $\alpha, \beta\in \mathbb{R}$ ($\alpha, \beta, \gamma \in \mathbb{R}$, respectively) such that the formula \eqref{first} (\eqref{first1}, respectively)   makes sense,
   let us denote by $\hat g[g, \bar g, \alpha, \beta]$ ($\hat g[g, \bar g, \tilde g,  \alpha, \beta, \gamma]$, respectively)   the  metric \eqref{first} (\eqref{first1}, respectively): 
   \begin{equation} \label{first}
\hat g[g, \bar g, \alpha, \beta]:=   \frac{ {\alpha \cdot g}/{(\det(g))^{{2}/{3}}}+ {\beta\cdot \bar g}/{(\det(\bar g))^{{2}/{3}}}}{\left(\det\left( {\alpha\cdot  g}/{(\det(g))^{{2}/{3}}}+ {\beta\cdot \bar g}/{(det(\bar g))^{{2}/{3}}}\right)\right)^2}
   \end{equation}
    \begin{equation} \label{first1}
\hat g[g, \bar g, \tilde g, \alpha, \beta, \gamma]:=   \frac{ {\alpha \cdot g}/{(\det(g))^{{2}/{3}}}+ {\beta\cdot \bar g}/{(\det(\bar g))^{{2}/{3}}} + {\gamma \cdot \tilde g}/{(\det(\tilde g))^{{2}/{3}}}}{\left(\det\left( {\alpha\cdot  g}/{(\det(g))^{{2}/{3}}}+ {\beta\cdot \bar g}/{(det(\bar g))^{{2}/{3}}}+ {\gamma \cdot\tilde  g}/{(\det(\tilde g))^{{2}/{3}}}\right)\right)^2}
   \end{equation}
   In these   formulas, $g$, $\bar g$,   and $\tilde g$ should be understood as $(2\times 2)$-matrices of metrics in a local coordinate system. In  \S\ref{2.1} and \S\ref{4.1},  we will explain the geometry and the hidden linear structure   behind this formula. In particular, it will be clear that the  formula is independent of the coordinate system (though one can check it by hand).
   Moreover,  if the metrics $g$ and $\bar g$  ($g$, $\bar g$,   and $\tilde g$, respectively) 
       are projectively equivalent, then $\hat g[g,\bar g, \alpha, \beta]$  ($ \hat g[g, \bar g, \tilde g, \alpha, \beta, \gamma]$, respectively) is also projectively equivalent to $g$. 
   Of course,  the metrics $\hat g[g, \bar g, \alpha, \beta]$  ($ \hat g[g, \bar g, \tilde g, \alpha, \beta, \gamma]$, respectively) 
 are  defined only for $\alpha, \beta, \gamma  \in \mathbb{R}$ such that  $\det\left( {\alpha\cdot  g}/{(\det(g))^{{2}/{3}}}+ {\beta\cdot \bar g}/{(\det(\bar g))^{{2}/{3}}}\right)\ne 0$ or respectively$\det\left( {\alpha\cdot  g}/{(\det(g))^{{2}/{3}}}+ {\beta\cdot \bar g}/{(\det(\bar g))^{{2}/{3}}}+ {\gamma \cdot \tilde g}/{(\det(\tilde g))^{{2}/{3}}}\right)\ne 0$.

  Denote by  $G(g, \bar g)$   ($G(g, \bar g, \tilde g)$, respectively) the following set of  metrics:  
  \begin{equation}\label{G} G(g, \bar g):=  
\left\{  \hat g [g, \bar g, \alpha, \beta]     \mid \textrm{
  $\alpha, \beta \in \mathbb{R}$ such that $\hat g[g, \bar g, \alpha, \beta]$ is defined 
  }
  \right\}.\end{equation}
  \begin{equation}\label{G1} G(g, \bar g, \tilde g):=  
\left\{  \hat g [g, \bar g, \tilde g,  \alpha, \beta, \gamma ]     \mid \textrm{
  $\alpha, \beta, \gamma  \in \mathbb{R}$ such that $\hat g[g, \bar g, \tilde g,  \alpha, \beta, \gamma]$ is defined 
  }
  \right\}.\end{equation}
 As we explained above,
  if the metrics $g$ and $ \bar g$ ($g,\, \bar g,$ and $ \tilde g,$ respectively)  are projectively equivalent, then  $ G(g, \bar g)$  ($G(g, \bar g, \tilde g)$, respectively) is a subset of their projective class.

 Now, every  metric   $g$ from Theorem \ref{main} always admits a  nontrivial projectively equivalent metric: as we explain in Appendix, for arbitrary data $X(x)$, $X_1(x)$,  $Y(y)$, $Y_1(y)$,  $h(z)$, $h_1(z)$,    \begin{itemize} 
 \item  the metric $g $  from  the Liouville Case  
 of Theorem  \ref{main} is projectively equivalent to the metric 

 \begin{eqnarray} \label{E1}   ds_{\bar g}   & =&  \left(\frac{1}{X(x)}- \frac{1}{Y(y)}\right)\left(\frac{X_1(x)}{X(x)}dx^2+   \frac{Y_1(y)}{Y(y)}dy^2\right),\end{eqnarray}

 \item  the metric $g $  from  the Complex-Liouville Case  of Theorem  \ref{main} is projectively equivalent to the metric 
  \begin{eqnarray}   ds_{\bar g}   & =&   
        \left(\frac{1}{\overline{h( z})} - \frac{1}{h(z)} \right)\left(\frac{\overline{h_1( z})}{\overline{h( z})}d\bar z^2 - \frac{h_1(z) }{h(z)}dz^2\right),\end{eqnarray} \item  the metric $g $  from  the Jordan-block-Case  of Theorem  \ref{main} is projectively equivalent to the metric 
 \begin{eqnarray} \label{E3}  ds_{\bar g}   & =&  -\frac{2(Y(y)+x)}{y^3}dxdy  + \frac{(Y(y)+x)^2}{y^4}dy^2. \end{eqnarray} 
\end{itemize} 
 Such a metric $\bar g$ projectively equivalent to $g$ will be called a \emph{canonical projectively equivalent metric}\footnote{The notion is not coordinate independent and has sense only if the metric has the form as in Theorem \ref{main}}.

 Moreover, the metric $g$ from the  case \ref{3d} of Theorem~\ref{main} admits one more metric projectively equivalent to it and essentially different from  the canonically projectively equivalent metric given by \eqref{E3}, namely $\tilde g$ given by  \begin{equation} \label{tilde}  
 ds_{\tilde g}^2 =  {\frac {{9\,{\it dx}}^{2}}{ \left( {y}^{2}+x \right) ^{2}
 \left( 3\,x-{y}^{2} \right) ^{6}}}-4\,{\frac {y
 \left( 9\,x+{y}^{2} \right) {\it dx}\,{\it dy}}{ \left( 3
\,x-{y}^{2} \right) ^{6} \left( {y}^{2}+x
 \right) ^{3}}}+{\frac {12x{{\it dy}}^{2}}{ \left( {y}^{
2}+x \right) ^{2} \left( 3\,x-{y}^{2} \right) ^{6}}}
\end{equation}

 \begin{Th} \label{main1} \label{main2}  
 The projective class of every metric $g$  from  cases \ref{1a}--\ref{3c} of  Theorem~\ref{main} 
 coincides with $G(g,\bar g)$, where $\bar g$ is the canonical projectively equivalent  metric. 
 The projective class of the metric $g$ from the case \ref{3d}  of  Theorem~\ref{main} 
 coincides with $G(g,\bar g, \tilde g)$, where $\bar g$ is the canonical projectively equivalent  metric, and $\tilde g$ is the metric given by \eqref{tilde}.  
\end{Th} 

We see that Theorem \ref{main} describes all   projective classes admitting an essential  projective vector field, and Theorem \ref{main2}  describes all metrics in these projective classes. 
Theorem~\ref{main3} below implies that all  these  metrics  actually have $dim(\mathfrak{p})=1$, because by  \cite{bryant}  a metric admitting two independent projective vector fields  admits a Killing vector field.

\begin{Th} \label{main3} 
None of the metrics from Theorem \ref{main} admits a nontrivial  Killing vector field.
\end{Th}

 Theorems \ref{main}, \ref{main2}, \ref{main3} give a complete solution of the  Problem 1 (of Lie) above.

\subsection{New ideas compared with [9]}
The theory of projective transformations  and projectively equivalent metrics has a long and fascinating history. The first non-trivial
examples of projectively equivalent metrics and projective transformations were discovered by Lagrange   \cite{Lagrange} and  Beltrami \cite{Beltrami}.   
Recently,  there has been a considerable growth in interest in projective
differential geometry, due to  new methods   that allow one  to solve interesting new and classical problems, see for example \cite{BK,einstein,hyperbolic,topology,fomenko,oshemkov}. 

The present  paper also solves an interesting    classical problem     explicitly stated by  Sophus Lie.  In a certain sense, this  paper  is a continuation of  \cite{bryant},  where another problem stated by Sophus Lie (Problem 2 from \S \ref{introduction}) was solved; let us recall the main idea of \cite{bryant} and  comment in brief on  new  (with respect to \cite{bryant}) ideas  of the present  paper. 

   Let $S$ be {\it a projective structure}  on  a surface, i.e., 
    an equivalence class of torsion-free affine connections with the same unparameterized geodesics (in \S \ref{schema} we will explain that  in a coordinate system  projective structures are parametrized by four functions).  Certain  projective structures contain
     the Levi-Civita connection of a metric, in this case we say that the metric is {\it compatible} with the projective structure, and the projective structure is {\it metrizable}.    

    Projective structures with  projective vector fields  are easy to classify: in particular, projective structures admitting two projective vector fields were essentially described  by Lie himself \cite{Lie} and Tresse \cite{Tresse}.   In order to solve Problem 2 of Lie, one needs  to understand which  projective structures admitting two projective vector fields are metrizable.

By an old and now well-understood observation of R. Liouville \cite{liouville}, metrics  compatible with  a projective structure 
can be found as nondegenerate solutions of an overdetermined system  of linear  partial differential equations, whose coefficients  are given by  the projective structure. 
There exists an algorithmic way   (sometimes called the prolongation-projection method)  to  understand whether an  overdetermined system  of linear    partial differential equations has a nontrivial  solution. In \cite{bryant}, the     algorithm   was applied and   all metrics whose projective structures  admit  two projective vector fields were  described.  Moreover,  recently the algorithmic way  to understand whether a given  projective structure is  metrizable  
 was essentially     simplified in \cite{duna}.

Unfortunately,   it is hard to apply this method  to find all metrics admitting one projective vector field. Indeed,   projective structures admitting one  projective 
vector field  depend explicitly on arbitrary functions
of one variable. The prolongation-projection  method (or the results of \cite{duna}) applied to  such projective structures results in a completely intractable nonlinear system of ODEs, yielding no insight.

In order to solve Problem 1 of Lie we used   another method. 
We used the fact that  the system of PDEs   defining compatible metrics is projectively invariant, hence
the solution space is invariant under the Lie derivative by a projective vector field. If a metric admits an essential projective vector field, the family of compatible metrics  is  at least two-dimensional.  The question of finding projective structures with a 2-dimensional family of compatible  Riemannian metrics was posed by  Beltrami \cite{Beltrami} and solved by  Dini \cite{Dini}.   The solution
depends explicitly on two arbitrary functions, each of one variable. However, we  are interested in all signatures: even if the original metric $g$ is Riemannian,
other metrics in the family need not be. In signature (+,--), there are more
solutions: in addition to a straightforward  analogue of the Riemannian solution, there
is a ``complex" form and a ``degenerate" case. Although this may have been known
to  Darboux and other authors,  the appendix of the paper provides a straightforward
and self-contained description of pairs of projectively equivalent metrics.

Returning to the main thread of the paper, the strategy now is to analyse the linear
action of the Lie derivative along a projective vector field on the space of solutions to the
equation for compatible metrics. In the nontrivial case, this action turns out to be invertible
with a 2-dimensional invariant subspace, see \S \ref{ii}. The form of the metrics corresponding to this
subspace is thus given by one of the three cases in  the appendix. There are also three
possibilities for the linear action: it may have a single 1-dimensional eigenspace, complex
conjugate eigenvalues or two real eigenvalues.
Consequently there are nine separate cases to consider, each of which gives rise to a
system of six PDEs for  the at most  two unknown functions (of one variable) in the metrics and 
two unknown components  of the projective vector field. The Leibniz
rule for the Lie derivative implies that the system is linear in the six first derivatives
of these functions. It turns out that the system is  solvable for these first derivatives (with
independent variables not appearing explicitly). This leads to their explicit integration
in each case.    The result is the explicit classification  of  metrizable projective structures  
 admitting  at least one projective vector  field and  compatible with  at least two nonproportional metrics,   i.e.,  essentially,  Theorem \ref{main}.

However, this is not yet an explicit classification of metrics with a nontrivial vector
field, as an explicit form has only been found for one  metric in the same projective class.
The classification is completed by Theorem \ref{main1}, in which all metrics projectively equivalent
to those in Theorem \ref{main} are described, and Theorem \ref{main3}, which shows that the dimension
of the space of projective vector fields is exactly 1. The proof of Theorems \ref{main1},  \ref{main3} is standard (though quite nontrivial technically). In order to prove Theorem \ref{main1}, 
  we apply the (adapted version of the)  prolongation-projection method  to show that no other solution exists. 
  In order to prove  Theorem \ref{main3}, we use a certain trick  known to Darboux and Eisenhart, see Section \ref{pr4} for more details.

\subsubsection*{Acknowledgments}
I  would like to thank  R.~Bryant, 
 M. Dunajski,  M.~Eastwood,     and   B.~Kruglikov for useful discussions. I am very grateful to G. Manno for checking the paper and finding a lot of misprints, and to N. Hitchin for grammatical and stylistic corrections.   I thank  
Deutsche Forschungsgemeinschaft
(Priority Program 1154 --- Global Differential Geometry and Research Training Group   1523 --- Quantum and Gravitational Fields),
  and FSU Jena for partial financial support, and 
 Cambridge  and Loughborough   Universities   for hospitality.

\section{ Schema   of the proof of Theorem 1} 

Roughly speaking, we reformulate our problem as 9 systems  of  PDEs  and solve them. In this section we will explain how we do it. More precisely, \begin{itemize} \item in \S\ref{2.1},   we review  the  theoretical   results  we will use.\item  In \S\ref{nondegenerate},   we prove two additional (relatively simple) results.  

\item In \S\ref{trick},  we explain the main trick that allowed  us  
   to reduce our problem    to  9 systems of PDEs which are relatively easy and can be solved explicitly, possibly after an appropriate coordinate change. We will also explain in what sense the systems are easy.   

\end{itemize}

In Section \ref{calculation},   we solve  these  9 systems of PDEs. 

\subsection{\bf General theory }  \label{2.1} \label{schema} 
  General theory   can be found in   \cite{bryant, duna, MT,dim2, Quantum, quantum,   dedicata} and in   more classical sources which in particular 
  are acknowledged in  \cite{bryant}.  The present paper should be viewed as a continuation of  \cite{bryant},   it could be useful  for  a reader  to have  \cite{bryant}   at hand  while reading 
    the present  paper.    

We will work on a small disc $D^2$ in local coordinates $(x,y)$.

\begin{Def} A \emph{ projective connection}  is a  second order ordinary differential equation of the form
\begin{equation}\label{eq1}
y'' = K_0(x,y) + K_1(x,y)\,y' + K_2(x,y)\,(y')^2 + K_3(x,y)\,(y')^3,
\end{equation}
where the functions~$K_i:D^2\to \mathbb{R}$. 

For any symmetric affine connection $\Gamma =(\Gamma^i_{jk}(x,y))$,
the projective connection \emph{ associated to}~$\Gamma$
is 
\begin{equation} \label{con}
y'' =-\Gamma^2_{11}+(\Gamma^1_{11}{-}2\Gamma^2_{12})\,y'
       -(\Gamma^2_{22}{-}2\Gamma^1_{12})(y')^2
         +\Gamma^1_{22}(y')^3.        
\end{equation}

We say that a  metric $g$ belongs to the 
\emph{ projective class of the projective connection \eqref{eq1}}, if  the projective connection \eqref{eq1}
 is  associated to the Levi-Civita connection of $g$.
\end{Def}

As has been known since the time of Beltrami~\cite{Beltrami},  
the projective connection associated to $\Gamma$ 
 carries  all the information about unparameterized geodesics of $\Gamma$. More precisely, for every solution $y(x)$ of~\eqref{con}, 
the curve $\bigl(x,y(x)\bigr)$ is, up to reparameterization,
a geodesic of~$\Gamma$. 
In particular, two metrics are projectively equivalent if and only if  they belong to the  projective class of the same projective connection. Therefore, according to the definition in  \S\ref{introduction},  
the   { projective class of $g$}  
is the projective class of the   projective connection associated to the Levi-Civita connection of $g$.

Let us reformulate (following  \cite{liouville,bryant}) the condition \\[.5ex]  {\em ``the metric $\  \  E(x,y)\,dx^2+2F(x,y)\,dx\,dy+G(x,y)\,dy^2$ belongs to the projective class of \eqref{eq1}}" \\[.5ex]  as a system of PDE on $E, F, G$.

Consider  the symmetric nondegenerate matrix
\begin{equation} \label{a}
a = \left(\begin{array}{cc}a_{11}&a_{12}\\a_{12}&a_{22}\end{array}\right)
:= \det(g)^{-2/3} \cdot  g
 = \frac1{(EG-F^2)^{2/3}}\left(\begin{array}{cc}E&F\\F&G\end{array}\right).
\end{equation}

\begin{Lemma}[\cite{liouville,bryant}]
\label{1}
The projective connection { associated to} the Levi-Civita connection 
 of the metric $g$ is~\eqref{eq1}
if and only if the entries of the matrix~$a= \det(g)^{-2/3} \cdot  g$ 
satisfy the \emph{linear} PDE system
\begin{equation}
\label{lin1}
\left.
\renewcommand{\arraystretch}{1.5}
\begin{array}{rcc}
{a_{11}}_x-\tfrac{2}{3}\,K_1\,a_{11} +2\,K_0\,a_{12}&=&0\\
{a_{11}}_y+2\,{a_{12}}_x
-\tfrac{4}{3}\,K_2\,a_{11}+\tfrac{2}{3}\,K_1\,a_{12}+2\,K_0\,a_{22}&=&0\\
2\,{a_{12}}_y+{a_{22}}_x
-2\,K_3\,a_{11}-\tfrac{2}{3}\,K_2\,a_{12}+\tfrac{4}{3}\,K_1\,a_{22}&=&0\\
{a_{22}}_y-2\,K_3\,a_{12}+\tfrac{2}{3}\,K_2\,a_{22} &=&0
\end{array} \right\}
\end{equation}
\end{Lemma}

In the formula \eqref{a}  above, $a$ should be understood as a section of \begin{equation} \label{sec} S_2D\otimes 
\left(\Lambda_{2} D\right)^{-\frac{4}{3}},\end{equation} 
where $\Lambda_2$ is the one-dimensional bundle of volume forms. 
Indeed, after a coordinate change the transformation rule  of an element of \eqref{sec} and of $\det(g)^{-2/3} \cdot  g$ coincide. 

In particular, the Lie derivative of $a= \det(g)^{-2/3} \cdot  g$ is well defined (as a mapping from the sections of $ S_2D\otimes 
\left(\Lambda_{2} D\right)^{-\frac{4}{3}}$ to itself), is independent on the coordinate  system,  and is given by the formula
\begin{equation} \label{liederivative} 
L_v a =L_v \left(\det(g)^{-2/3} \cdot  g\right)  = \det(g)^{-2/3} \cdot L_v  g - \frac{2}{3} \det(g)^{-2/3} \trace_g(L_v g )\cdot g  ,   \end{equation}  
where $\trace_g(L_v g ) := \sum_{i,j}\left(L_v g\right)_{ij}g^{ij }.$

\begin{Rem} \label{rem3} 
The formula \eqref{a} is invertible:  
$a = g/\det(g)^{2/3}$ if and only iff 
 $g= a/\det(a)^2$.   The mapping $a\mapsto  a/\det(a)^2$ can be viewed as a mapping from  $S_2D \otimes\left(\Lambda_{2}D\right)^{-\frac{4}{3}}$ to  $S_2D$, which  is defined for nondegenerate 
  points of $S_2D \otimes\left(\Lambda_{2}D\right)^{-\frac{4}{3}}$ only, and sends them into  nondegenerate 
  points of $S_2D $.  In particular, if  a nondegenerate 
  $a$ is a solution of \eqref{lin1}, then the metric  $g= a/\det(a)^2$ belongs  to the  
   projective class of \eqref{eq1}.
  \end{Rem} 

The system \eqref{lin1} has the following  nice properties, which will be used later: 
\begin{itemize} 
\item It is linear and of finite type.  In particular, its space of solutions  (that will be denoted by $\mathcal{A}$) is a finite-dimensional ($\dim(\mathcal{A})\le 6, $ \ \cite{liouville} )  vector space.  
\item  Moreover, if $dim(\mathcal{A})\ge  4$, then every metric from the projective class admits a Killing vector field \cite{koenigs}.   
\item The system \eqref{lin1}  depends on the projective connection only and is therefore projectively invariant. In particular, for every projective vector field $v$ and for every solution $a\in \mathcal{A}$ we have $L_v a \in \mathcal{A}$.   Thus, $L_v$ is a (linear)  mapping from  $\mathcal{A}$ to  itself.   
\end{itemize} 

We will also use the following two statements: the first is due to Knebelman  \cite{knebelman} (another proof can also be found in  \cite{bryant,kruglikov,MT,dedicata,beltrami_short, BMF},   one more 
proof easily  follows from the theory of  invariant  operators, see for example  \cite{baston}).  The second  is a combination of the formula \eqref{a} and  the connection between projectively equivalent metrics and integrable systems 
due to \cite{dim2,MT},    see also Darboux~\cite[\S608]{Darboux}, see also of \cite[\S 2.4]{bryant}.

\begin{itemize}
\item If a metric $g$ admits a Killing vector field, 
then every metric projectively equivalent to $g$ also admits a Killing vector field. 

\item  $a$ is a solution of the system \eqref{lin1} corresponding to  (the projective connection associated to the Levi-Civita connection of) a metric $g$,  if and only if 
  the function \begin{equation} 
I:TD^2 \to \mathbb{R}\ , \ \  I(\xi)= \det(g)^{2/3} \cdot  \sum_{i,j} a_{ij}\xi^i \xi^j \label{integral} \end{equation} is an integral for the geodesic flow of $g$.  
\end{itemize} 

\begin{Rem} { \ }

\begin{enumerate}

 \item Tensor products with powers of  $\left(\Lambda_{2} D\right)^{\frac{1}{3}}$ appear naturally in the theory of projectively equivalent metrics and projective connections, see \cite{eastwood1}. 

 \item A multidimensional generalization of the formula   $a := \det(g)^{-2/3} \, g$ and  of Lemma \ref{1}  can be found in \cite{eastwood}, see also \cite{Aminova, benenti, mikes, sinjukov}.

 \item The formula \eqref{liederivative} appears naturally in the  investigation of projective transformations of surfaces, see  \cite{obata,sol, CMH}, and can be generalized to all dimensions, see \cite{archive, topalov}.   \end{enumerate}
\end{Rem} 

\subsection{  Every nontrivial solution $a$ of the system (12) is nondegenerate at almost every point, and $L_v:\mathcal{A}\to \mathcal{A}$ is nondegenerate. }  \label{ii} 

Within this paragraph  we assume that   the restriction of  $g$  
 to every open  neighborhood $U\subseteq D^2$  admits no Killing vector field. We 
denote by $\mathcal{A}$ the space of  solutions of the system \eqref{lin1}   corresponding to the projective connection associated to the Levi-Civita connection of  $g$. 

\begin{Lemma} \label{nond}  Assume  $a\in \mathcal{A}$ is not identically zero. Then, 
 the set of the points where $a$ is degenerate  is nowhere dense (in the topological sense, i.e.,  the complement to  this set is open and everywhere dense.). 
\end{Lemma} 
\noindent  {\bf   Proof.}  The set of  points where $a$ is degenerate  is evidently a closed set.  Assume  there exists a neighborhood $U$ such that  $a$ is degenerate at every point of $U$. In a local coordinate system $(x^1,x^2)$ in the neighborhood $U$  $a$ is given by a   symmetric $(2\times 2)-$matrix $a= \begin{pmatrix} a_{11}  &  a_{12} \\ a_{12} & a_{22}\end{pmatrix}$. If the kernel of $a$ is two-dimensional at every point of  a certain 
 neighborhood  $U\subseteq D^2$, then the restriction of $a$ to $U$ is identically zero.   
 Since the PDE system \eqref{lin1} is linear and is  of  finite type,    
 $a\equiv 0$  on  the whole  $D^2$. Then, the set of the points where the kernel of $U$ is precisely one-dimensional is everywhere dense in $U$.  
  Without loss of generality, passing to a smaller neighborhood if necessary,   we may assume that the kernel of $a$ is precisely  
  one-dimensional at every point of   $U$. Take a local coordinate system $(x,y)$
   on an open subset  $U'\subseteq U$  such that the kernel of $a$   is the linear hull of 
    $\frac{\partial }{\partial y}$. Then,  in this coordinate system 
     the matrix $a$ has the form $a= \begin{pmatrix} \alpha(x,y)  &  0 \\ 0 & 0\end{pmatrix}$, 
     where $\alpha$ vanishes at no point of $U'$.  

  Then, the  integral \eqref{integral}  of the geodesic flow of $g$ is equal to $\det(g)^{2/3} \cdot  \alpha \left(\xi^1\right)^2$. Then, the function $I_{lin}:= \sqrt{\det(g)^{2/3} \cdot  |\alpha|} \xi^1$ 
  is also an integral. Since the integral $I_{lin}$  is linear in velocities, 
  the metric $g_{|U'}$  admits a Killing vector field. The contradiction proves Lemma~\ref{nond}.  \qed

\label{nondegenerate} 
\begin{Lemma} \label{kil} For every projective vector field $v$,  
the mapping $L_v:\mathcal{A}\to \mathcal{A}$ is nondegenerate. 
\end{Lemma}  

\begin{Rem} See \S\ref{2.1}  for interpretation $L_v$ as a linear mapping from $\mathcal{A}$ to $\mathcal{A}$. 
\end{Rem} 
\noindent{\bf Proof of Lemma \ref{kil}.} Assume  there exists an 
nontrivial $\bar a\in \mathcal{A}$ such that 
$L_v\bar a = 0$. In a neighborhood of the point such that $v\ne 0$ take a coordinate system $(x,y)$  such that $v= \frac{\partial }{\partial x}$. Then, the components of 
$ L_v\bar a$  are the $x-$ derivatives of the components of $\bar a$, and the condition $L_v\bar a = 0$ implies that the components of $\bar a$ are independent of $x$. Then, the components of the metric $\bar g:=  {\bar a}/{\left(\det(\bar a)\right)^2}$, which is defined almost everywhere by Lemma \ref{nond}, are independent of $x$. Then, $v$ is a Killing vector field for $\bar g$.  Since, as we explained in \S\ref{2.1}, see Remark \ref{rem3} there, 
the metric $g$ is projectively equivalent to $\bar g$,    
then, by result of Knebelman \cite{knebelman} we recalled in \S\ref{2.1},  the metric $g$ also  admits a Killing vector field in a neighborhood of almost  every point. The contradiction proves  Lemma \ref{kil}. \qed

\subsection{\bf How to reduce Theorem~1  to $9$ Frobenius  systems of PDEs  } \label{trick} 

 Recall that a PDE-system of the first order is \emph{ Frobenius}, if the derivatives of all unknown functions are explicitly  given as functions of the unknown functions. Frobenius systems are easy to handle: there exists an algorithmic way to reduce them to ODEs. In our case, the Frobenius systems are  simple enough so that we  could  explicitly solve them.  Note that the most straightforward way   to reformulate the problem as a system of PDEs, i.e., if we write down the  conditions that a vector field $\frac{\partial }{\partial x}$ is projective with respect to $g$ as a system of 4 PDEs in 3  
 unknown components of the metric,  leads to a much more complicated system of PDEs which is  impossible  ( = we did not find a way to do it) 
 to solve by standard methods. In fact, the system is only slightly overdetermined (4 equations on 3 unknowns), and the standard prolongation-projection method  will require too many  (more than 20) operations of prolongation and prolongation-projections.

The reduction of Theorem~\ref{main} to 
$9$ Frobenius  systems of PDEs   is based on the description of projectively equivalent metrics in the appendix, and on  the following two trivial  statements from  linear algebra: 

\begin{itemize}

\item   For every  nondegenerate 
linear mapping $L:\mathbb{R}^2  \to  \mathbb{R}^2$ there exists a basis in $\mathbb{R}^2$ such that  for the appropriate $\const\in \mathbb{R}$ 
the matrix of $\const \cdot L$ is given  by 
\begin{equation}  (a) \begin{pmatrix} 1  & 1 \\ & 1 \end{pmatrix}, \  \   (b)
     \begin{pmatrix} \lambda & -1  \\ 1   &  \lambda\end{pmatrix} , \textrm{  or   } (c)           \begin{pmatrix} \lambda & \\ & 1\end{pmatrix}. \label{mat}\end{equation} Moreover in the matrix $(c)$ we can assume  $\lambda \in (-\infty, -1] \cup [1, \infty )$.  
\item Every nondegenerate  linear mapping $L:\mathbb{R}^3  \to  \mathbb{R}^3$ 
has a two-dimensional invariant subspace such that the matrix of the 
 restriction of $\const\cdot L$ to this subspace   is one of the matrices \eqref{mat}  in a certain basis.   
\end{itemize}
Let us  explain how we reduced  Theorem~\ref{main} to   (solving of) 
nine  Frobenius  systems of PDEs.  

Suppose the metric $g$ has an essential projective vector field $v$ and admits no Killing vector field. Consider the projective connection associated to the Levi-Civita connection  of the metric,   and the  space $\mathcal{A}$ of  solutions of  
 \eqref{lin1}.  
  Since the system \eqref{lin1} is projectively  invariant, for every  $a\in \mathcal{A}$     the  Lie derivative $L_v a$  is also a solution. Thus, $L_v$ can be viewed as   a linear mapping $L_v:\mathcal{A}\to \mathcal{A}$.

  {\bf The case  $dim(\mathcal{A})= 1$} is not interesting for us. Indeed, in this case, 
     all metrics  projectively equivalent to $g$ have the form $\const \cdot g$, which in particular implies  that all projective vector fields are  infinitesimal homotheties,
     and in Theorem~\ref{main} we excluded such metrics. 

    {\bf  The case  $dim(\mathcal{A})\ge 4$ } is also not interesting for us. Indeed, in this case,  
     as we   recalled  in \S\ref{2.1},   the metric  $g$ 
  admits a Killing vector field.

  { \bf  If   $dim(\mathcal{A})=2 \textrm{ or } 3$,}  then, as we explained above,  
    $\mathcal{A}$ has a two-dimensional invariant subspace $\hat{\mathcal{A}}$ such that the restriction of $L_v$ to $\hat {\mathcal{A}}$ is given by one of the matrices \eqref{mat}. 
  If   $\{a, \bar a\}\in \hat{\mathcal{A}}$  is the basis such that 
     $L_v$ is given by, say,  the matrix  (b)  from \eqref{mat}, we have  (the matrices (a) and (c) 
      will be treated in \S\ref{pa} and  \S\ref{pc}, respectively) 
  \begin{equation}\left. \begin{array}{ccc} L_v a  & = & \lambda  a - \bar a   \\  L_v\bar a  & = &   a + \lambda \bar a.\end{array}  \right\}\label{cond2}  \end{equation}

  By Lemma \ref{nond} from \S \ref{ii}, without loss of generality 
   we can assume  that the   the matrices of $a, \bar a$ are nondegenerate, since they are so at  almost every point.  Then, $a$ and $\bar a $ come from two certain  metrics by formula \eqref{a}, see Remark \ref{rem3}.  By Lemma~\ref{1}, the metrics are projectively equivalent to $g$; without loss of generality we can think that the metric corresponding to $a$ is the initial metric $g$. We will call $\bar g$ the metric corresponding to $\bar a$, so that 
   $$
   a = \det(g)^{-2/3} \cdot  g\ , \  \   \ \bar a= \det(\bar g)^{-2/3} \cdot  \bar g.
   $$
  Then, in view of \eqref{liederivative},  the condition \eqref{cond2} reads 
  \begin{equation}\left. 
  \begin{array}{ccccccc} L_v  g   & = &\frac{2}{3} \trace_g(L_v g) g &+  & \lambda  g  & - &   \left(\frac{\det(g)}{\det(\bar g)}\right)^{2/3}  \bar  g    \\ L_v  \bar g   & = &\frac{2}{3} \trace_{\bar g}(L_v \bar g)\bar  g &+&   \left(\frac{\det(\bar g)}{\det(g)}\right)^{2/3}    g & + &  \lambda  \bar g.    \end{array}\right\}  \label{cond3}  \end{equation}

 On the other hand, by  Theorem A from  the appendix, there exists a coordinate system $(x,y)$ such that the metrics $g$ and $\bar g$ are given by one of the model forms. Substituting the model metrics $g, \bar g$ from  the appendix, we obtain the system    of $6= 3+ 3$ PDEs \footnote{each  of two equations of    \eqref{cond3} is an equation on a symmetric $(2\times 2)-$matrix, i.e., is equivalent to three scalar equations}  of the first order on the data of the metrics and on the unknown projective vector field $v$.  

 Let us now count the number of first  derivatives of the unknown functions  in this system of  6 PDEs.
 In every model case, the data of  metrics  $g$ and $\bar g$, i.e., $X$ and $Y$ in the Liouville case, $h$ in the Complex-Liouville Case, and $Y$ in the Jordan-block Case,  have at most two first derivatives.  Together with four derivatives of the components of the vector field $v$, it gives us at most 6 first derivatives of the unknown functions. 

 Thus, in  the system  \eqref{cond3} the number of highest ( = first) derivatives is not greater than the number of the equations. It appears that in all cases it is possible to solve\footnote{since the system
 \eqref{cond3} and its analog for  matrices (a), (c) of  \eqref{mat}  is  linear in the derivatives, it is an exercise in linear algebra} the system with respect to the first derivatives, i.e., to bring the systems into  the Frobenius form,
  and then to solve it using standard methods.  

 We see that we have three choices for the matrix from \eqref{mat}, and three choices for the model metrics $g, \bar g$. Thus, we have $3\times 3=9$ Frobenius systems
 to solve. We  will   do so in Section  \ref{twodiminvarinat}.

 \begin{Rem}  In a certain sense, some systems from these nine are closely related, and can be obtained one from another by a kind of complexification. Indeed, as Remark A from  the appendix shows, the Complex-Liouville case could be understand as the complexification of the Liouville case. Moreover, over the field of complex numbers, the matrix (a) from \eqref{mat} has the same type as the matrix (b): they both have two different eigenvalues. One can indeed formalize these arguments  and reduce the number of systems to solve to four. But the nine systems are so simple that it is shorter to solve them than to explain how to make   a  solution of one using a solution of another. \end{Rem}

     \label{2.2}

\section{Calculations related to proof of Theorem~1} \label{calculation}  Within the whole section we assume that 

\begin{itemize} \item  $D^2$ is  a smooth disc with a  (Riemannian or pseudo-Riemannian) metric $g$ and  coordinates $(x,y)$, 
\item   the smooth vector field  $v$ is projective with respect to the metric $g$, 
\item the restriction of the metric $ g$  to any  open  subset $U\subset D^2 $ 
admits no Killing vector field.
\end{itemize}  

Within the whole section we work in the coordinates   $(x,y)$;  $g$ will always denote the metric we work with, and   
  $v=(v_1, v_2)$ its  projective vector field. As in \S  \ref{introduction}, we reserve the notation    $\varepsilon$ for $ \pm 1$. 

We consider the projective connection \eqref{eq1} associated to the Levi-Civita connection of the 
metric $g$, and       denote by $\mathcal{A}$ the space of  solutions of the equation \eqref{lin1}. We 
assume $dim(\mathcal{A})= 2 \textrm{  or } 3$ (see   \S \ref{2.1} for an explanation of why we can do this).

    \label{twodiminvarinat}

Let   $\hat{\mathcal{A}}\subseteq \mathcal{A} $ be a two-dimensional  subspace 
  invariant with respect to the  Lie derivative: $L_v(a)\in \hat{\mathcal{A}}$ for every $a \in \hat{\mathcal{A}}$ (we explained its existence in \S\ref{trick}).

Then, in view of Lemma \ref{kil}   and  
 after the multiplication of $v$ by an appropriate nonzero constant,  in a certain basis $\{a , \bar a\}$ of $\hat{\mathcal{A}}$ the matrix of $L_v$ is as in \eqref{mat}.

  We have three possibilities for the matrices in \eqref{mat}, we will carefully consider them in \S\S\ref{pa}, \ref{pb}, \ref{pc}.

  \subsection{ The matrix  of $L_v$ is as (a) in (16) }   \label{pa} 

Assume that in the basis $\{a  , \bar a\}$ the matrix of $L_v:\hat{\mathcal{A}}\to \hat{\mathcal{A}}$ is given by $\begin{pmatrix} 1 & 1 \\ & 1 \end{pmatrix}.$

Without loss of generality, in view of Lemma~\ref{nond} and Remark \ref{rem3}, we can assume that $
a = \det(g)^{-2/3} \cdot  g\ , \  \   \ \bar a= \det(\bar g)^{-2/3} \cdot  \bar g$ for certain metrics $g, \bar g$ from the projective class of \eqref{eq1}. 
Then, as we explained in  \S\ref{trick}, the condition 
\begin{equation} L_v\begin{pmatrix} a \\ \bar a\end{pmatrix} =  \begin{pmatrix} 1 & 1 \\ & 1   \end{pmatrix} \begin{pmatrix} a \\ \bar a\end{pmatrix} \textrm{ \  or, equivalently, \ }  \left\{ \begin{array}{ccccc} L_v a  & = & a &+&  \bar a   \\  L_v \bar a  & = &   &&  \bar a.\end{array} \right. \label{17} \end{equation}
is equivalent to  the following condition: 
\begin{equation}\label{en1t}\left.\begin{array}{ccccr} L_v g  - \frac{2}{3} \trace_g (L_v g) \cdot g & =  &g  &+ & \left(\frac{ \det(g)}{\det(\bar g ) }\right)^{2/3} \bar g     \\
   L_v \bar g  - \frac{2}{3} \trace_{\bar g} (L_v \bar g) \cdot \bar g & = &   && \bar g.    \end{array}\right\}\end{equation}

  As  explained in  the appendix, 
   in a neighborhood of almost every point the metrics $g$ and $\bar g$ have one of three normal forms. We will carefully consider all three cases. 

\label{J}

\subsubsection{ Liouville Case}    \label{La}
   Assume they have the Liouville form, i.e., 
  \begin{equation} \label{liu} 
  ds^2_g=(X-Y)(dx^2 + \varepsilon     dy^2), \  \  \    ds^2_{\bar g}=\left(\frac{1}{Y}-\frac{1}{X}\right)\left(\frac{dx^2}{X} +  \varepsilon   \frac{dy^2}{Y}\right). 
  \end{equation}

    After some calculations we obtain that the Lie derivatives of $g$ and $\bar g$ are  given by   the matrices  
    $$ \left(\begin{array}{cc}   {X'} {v_1} +2\, X \frac{\partial  v_1}{\partial x}    -2 Y \frac{\partial v_1}{\partial x}    
  -  {Y'}  v_2   & 
  \left( {\frac {\partial  v_1}{\partial y{{}}}}{  }  + \varepsilon {\frac {\partial v_2 }{\partial x{{}}}}{   }   \right)  \left( X  -Y
   \right) \\  \left( {\frac {\partial   v_1}{\partial y{{}}}}{ } +\varepsilon {\frac {\partial  v_2}{\partial x{{}}}}{  }   \right)  \left( X \ -Y
   \right)&   \varepsilon\left( {{X'}}    {   v_1}  - {{Y'}{{{}}}}
   {   v_2}  +2\,  X {\frac {\partial {   v_2}}{\partial y{{}}}}
     -2\,Y
 {\frac {\partial  v_2}{\partial y{{}}}}\right)    \end{array}  \right),$$ 
$$
\left( \begin {array}{cc} {\frac {YX'v_1 (2 Y-  X) +
2XY\frac{\partial v_1 }{\partial x}\left(X-    Y\right)
-  Y' X^2   v_2  }{  Y^{2} X^{3}}}
&
\frac { \left( X  -Y  \right)  \left(   
\frac {\partial v_1 }{\partial y}  Y  + \varepsilon \frac{\partial v_2 }{\partial x}    X\right) }{  Y^{2}  X^{2}}
\\ \frac{ \left( X    -Y
   \right)  \left(   \frac{\partial v_1 }{
\partial y}     Y
 +  \varepsilon \frac{\partial  v_2}{\partial x}X   \right)}{  Y^2 X^2}&
 \varepsilon\frac{X'   Y^2v_1+  XY'v_2 (Y   
  -2 X  )+2     X Y (X-Y)   \frac {\partial v_2}{\partial y} }{ Y^{3}
  X^2}
\end {array} \right)  
$$
  and   the system \eqref{en1t}   is equivalent to the following  system of 6 PDEs in the unknown functions $v_1(x, y)$, $v_2(x, y)$,  $X(x)$, and $Y(y)$. 
\begin{equation} \left.  \begin{array}{ccc} 
  \frac{Y'v_2}{3}  -\frac{ X'v_1}{3} +\frac{2}{3}  \frac {\partial v_1}{\partial x}   (X  - Y)   
  -\frac{4}{3}  \frac{\partial v_2}{\partial y} (
  X -Y)  &=& (Y +1)(X -Y)
     \\
  (X-Y) \left(\frac{\partial v_1}{\partial y} + \varepsilon  \frac{\partial v_2}{\partial x}\right)&=&0\\
  \frac{X'v_1 }{3}\,    -\frac{Y'   v_2}{3}  -\frac{2}{3}  \frac{\partial v_2 }{\partial y}\left(  X -   Y\right) +\frac{4}{3}  
\frac{\partial v_1}{\partial x}   (X  -  Y) &=&  (1+ X) (Y-X) 
\\
  Y X' v_1 - 2\,Y \frac{\partial v_1}{\partial x}\left(   X-    Y\right)
  -v_2Y'(3  X    -2    Y)  +4\,  Y  \frac{\partial v_2}{\partial y}\left( X -  Y\right)
  &=& 3Y  (Y -  X)  \\
 \left( X  -Y  \right) 
\left(  \frac{\partial v_1}{\partial y} Y   + \varepsilon \frac {\partial v_2 }{\partial x}   X  
 \right)&=&0 
  \\
X'   v_1(3Y    -2  X)- 
Y'  X  {v_2} 
-2  X \frac{\partial v_2}{\partial y}\left( X-  Y \right)
+4 \frac{\partial v_1}{\partial x}   X(X
-Y) &=& 3X(Y-X) 
\end{array} \right\} \label{all} \end{equation}
We see that (in view of the nondegeneracy of the metric $(X-Y)(dx^2 + \varepsilon dy^2)$) the second  and fifth equations  of \eqref{all}  imply that $v_1$ depends on  the variable $x$ only, $v_2$ depends on  the variable $y$ only.    Then, all unknown functions in the  system  \eqref{all} are functions of one variable only, so  the system (\ref{all}) is actually a system of ODEs  (of first order). We see that it is 
   linear in the derivatives.
  Solving it for the derivatives of the unknown functions $X(x), Y(y), v_1(x), v_2(y)$, we obtain that   \eqref{all} is equivalent to the following  4 ODEs: 
  \begin{equation} \label{ode}
    v_1'  =-\frac{X}{2}-\frac{3}{2}\   , \  \  \ v'_2  =-\frac{3}{2}-\frac{Y}{2}   \  , \ \ \   Y' v_2 = -Y ^{2}\ , \ \ 
 \ X'v_1  =-X^{2}.  \end{equation}

These equations can already be solved; since the solution  is quite complicated and is given in terms of  Lambert functions,  instead of solving the system we change the coordinates  (possibly passing to a smaller neighborhood)  such that in the new coordinates the metrics $g$ and $\bar g$ and the vector field $v$ are given by  elementary functions.

Since by assumption the metric $g$   admits no Killing vector field, the functions $X,Y$
 are not constant in every neighborhood,  which in particular implies  that for almost every point  the functions $v_1, v_2$ are not zero in a neighborhood of the point. In such a neighborhood consider the 
 coordinate change   \begin{equation} (x,y)= \bigl(x(x_{old}), y(y_{old})\bigr) \ \textrm{given by} \ \   dx =\frac{1}{v_1} dx_{old} \   \   dy =\frac{1}{v_2} dy_{old}.\label{dx2} \end{equation} 
After this coordinate change the ``old" equation  $ X'v_1  =-X^{2}$ reads 
$\dot X = -X^2$, where $\dot X := \frac{d X }{d x}$, $\dot Y := \frac{d Y }{d y}$, $\dot v_1 := \frac{d v_1 }{d x}$, $\dot v_2 := \frac{d v_2 }{d y}$ are  derivatives with respect to the new coordinates.  The equation can be solved,  its nonconstant solution is 
$X(x) = \frac{1}{x+ c}$. Since the  formula \eqref{dx2} defines the coordinates 
 up to addition of arbitrary  
 constants,  without loss of generality we assume $c=0$, so that $X= \frac{1}{x}$. 
 Similarly, in the new coordinates the equation $ v_1'  =-\frac{X}{2}-\frac{3}{2}$ reads 
 $
 2\dot v_1 = -({X}+{3})v_1.   
 $    After substitution 
 $X= \frac{1}{x}$ we obtain 
 $
 2\dot v_1 = -\left(\frac{1}{x}+{3}\right)v_1,  
 $
 which can be easily solved, the solution is 
 $(v_1(x))^2 = \frac{ C_1}{ x} \cdot e^{-3x}  $. 
 Similarly, in the new coordinates the functions $v_2, Y$ are given by 
 $(v_2(y))^2 =  \frac{C_2}{y} \cdot e^{-3y} \   \  , \  \   \   \  Y(y)= \frac{1}{y}.$

 Thus, the metrics $g$  and $\bar g$ and the projective  vector field $v$ are  given by 
 \begin{equation*}  \begin{array}{ccccc} ds^2_{g}  & =  &
 (X-Y)(dx_{old}^2 \pm dy_{old}^2 ) & =& \left(\frac{1}{x}- \frac{1}{y} \right)\left( \frac{C_1}{x} e^{-3x}   dx^2 +  \frac{ C_2}{y} e^{-3y}   dy^2\right), \\ 
 ds^2_{\bar g}& = &  \left(\frac{1}{Y}-\frac{1}{X}\right)\left(\frac{dx_{old}^2}{X} \pm \frac{ dy_{old}^2}{Y} \right)    &=& \left(y- x \right)\left( C_1 e^{-3x}   dx^2 + C_2 e^{-3y}   dy^2\right), \\ v& =& v_1\frac{\partial }{\partial x_{old}} +   v_2\frac{\partial }{\partial y_{old}} & = &  \frac{\partial }{\partial x} +   \frac{\partial }{\partial y}. \end{array}
\end{equation*}

 We see that the metric $g$ and the vector field $v$  are as in case \ref{1a} of Theorem~\ref{main}. 

  \subsubsection{Complex-Liouville case } \label{CLa}

   Assume the metric $g$ and $\bar g$  have the Complex-Liouville form from Theorem A of   the appendix, i.e. 
  \begin{equation}\label{gcomplex}
 \begin{array}{ccc}  ds^2_g& = & 2  \Im(h) dxdy,\\  
   ds^2_{\bar g} & = & -\left(\frac{\Im(h)}{\Im(h)^2 +\Re(h)^2}\right)^2dx^2 +2\frac{\Re(h) \Im(h)}{   (\Im(h)^2 +\Re(h)^2)^2} dx dy   +  \left(\frac{\Im(h)}{\Im(h)^2 +\Re(h)^2}\right)^2dy^2.\end{array}  
    \end{equation} 

    \begin{Rem} 
    It could be helpful for understanding to know the complex version of the formulas \eqref{gcomplex}: it is 
   \begin{equation} \label{complexification} \left. \begin{array}{ccc}   ds^2_g & =  &  -\tfrac{1}{4}( \overline{h( z}) - h(z) )\left(d\bar z^2 -   dz^2\right)\\   
     ds^2_{\bar g} & =  &  -\tfrac{1}{4}\left(\frac{1}{\overline{h( z})} - \frac{1}{h(z)} \right)\left(\frac{d\bar z^2}{\overline{ h(z})} - \frac{ dz^2}{h(z)}\right),\end{array}\right\} \end{equation} 
     where $\bar z$ denotes the   complex-conjugate to $z$,   $\overline{h( z}) $ denotes  the   complex-conjugate to $h(z)$, and 
     $\bar g$  does not mean complex-conjugate to $g$, see Remark 1 from  Appendix. 

     We see that  the formula above is in a certain sense a complexification of \eqref{liu}, the role of $X(x)$ played by $h(z)$ and the role of $Y(y)$ by $h(\bar z)$.  We will see later, in all paragraphs related to the Complex-Liouville case, that all equations  related to the Complex-Liouville case could be viewed as a complexification of the corresponding equations from the Liouville case. Actually, one can show it advance, and avoid the calculation, but it appears that it is shorter to do the calculations than to explain why they could be avoided.   
    \end{Rem} 

    Arguing as in the previous paragraphs, we obtain that the 
    the conditions  \eqref{en1t}  are equivalent to a system of linear  6 PDEs of the first order. 
  Solving this system with respect to the first derivatives, and using the Cauchy-Riemann conditions for the holomorphic function $h$, we obtain  that the system is equivalent to the system 
  \begin{equation} \label{pp1} 
 \left.  \begin{array}{rcc} 
   v^2_y = {v^1_x}  &=&-\frac{\Re(h)}{2} -\frac{3}{2}\\ 
-v^1_y  =v^2_x  &=&-\frac{\Im(h)}{2}\\ 
 \Re(h)_x= 
 \Im(h)_y
  & =& \frac{ v^1    \left(\Im^2(h) -      \Re^2 (h)\right) -2 v^2 \Re(h) \Im(h) }{ \left(  v^2   \right)^2+ \left(v^1\right)^{2}}  \\
  -\Re(h)_y= { \Im(h)}_x & =&  -\frac{ v^2  \left( \Im^2(h)-  \Re^2(h)\right)
 + 2v^1 \Re(h)  \Im(h)}{ \left( v^2\right)^{2}+ \left( v^1   \right) ^{2}}
\end{array} \right\} \end{equation}

From the first two equations of \eqref{pp1}  we see that the function $V:= v^1 + i\cdot v^2$ is a holomorphic function of the variable $z:= x +i\cdot y$. It is easy to check that the last  two equations of \eqref{pp1} are equivalent 
to \begin{equation}  Vh_z= -h^2 \label{ppp1}, \end{equation} and the first two equations of   \eqref{pp1} are equivalent to \begin{equation}V_z= -\frac{h}{2} - \frac{3}{2}\label{pppp1} \end{equation}
 (where  $V_z$, $h_z$ are  the derivative of $V$ and  $h$ with respect to $z$). 
We see that  the equations (\ref{pp1} -- \ref{pppp1}) are direct analogues of  \eqref{ode}. 

After the holomorphic  coordinate change   
\begin{equation} dz_{\textrm{new}} =\frac{1}{V} dz_{\textrm{old}} \label{dz1n},  \end{equation}   
the equation  (\ref{ppp1}) is  $h_{z_{\textrm{new}}}=  -h^2$  implying
\begin{equation}h(z_{\textrm{new}})  =  \frac{1}{z_{\textrm{new}} + \const}\, \label{dz1bis} .\end{equation}
  Since the     formulas (\ref{dz1n}) define $z_{\textrm{new}}$ 
  up to addition of a complex  constant, we can (and will)  assume without loss of generality 
   that $\const=0$. In this new coordinate the vector field $v$ is $\frac{\partial }{\partial z} + \frac{\partial }{\partial \bar z}=\frac{\partial }{\partial x} $.

Now consider the equation (\ref{pppp1}). 
In the new coordinates it reads $$2{ V_{z_{\textrm{new}}}}{ }  = -
    \left( \frac{1}{z_{\textrm{new}}}  + {3}\right)V.$$
     Solving it, we obtain \begin{equation}\label{dz21}  {V}^2  
     ={\frac {{ C}  {e^{-3z_{\textrm{new}}}}}{ z_{{\textrm{new}}}
  }}.
\end{equation}  
Finally, substituting the  coordinate change and the solutions  (\ref{dz1bis} -- \ref{dz21}) in the metrics  we obtain that, after the appropriate scaling,  
the metrics $g$ and $\bar g$  have  the form 
 \begin{equation*}    \begin{array}{ccccl} ds^2_g & = &  2  {\Im\left({h}\right)}dx_{\textrm{old}} dy_{\textrm{old}}
 & = &   \frac{1}{4}\cdot{\left(h(\bar z)-h(z)\right)}  
\left(d\bar z_{old}^2  -  dz_{old}^2\right) \\ &  = &  && \frac{1}{4} \cdot \left(\frac{1}{\bar z} -\frac{1}{z}\right) 
\left(   \bar C  \, e^{-3 \bar z} \,  \frac{ d \bar z^2} {\bar z}
-    C  e^{-3 \,z}\, \frac{ d  z^2}{z} \right),\\
 ds^2_{\bar g} & = & &&  \frac{1}{4} \left(z- \bar z\right)\left(\bar C  \, e^{-3 \bar z}    d \bar z^2  -   C  e^{-3 \,z}\,  d  z^2\right), \end{array}\end{equation*} 
and the projective vector field $v $ is $\frac{\partial }{\partial z  } + \frac{\partial }{\partial \bar z }= \frac{\partial }{\partial x}$.

 We see that the metric and the projective vector field  $v$ 
  are as in case \ref{2a} of Theorem~\ref{main}. 

 \subsubsection{ The Jordan-block case} \label{JBa} 

  Let the metrics $g$ and $\bar g$ be given by the formulas from  Remark 2 of  the appendix:   \begin{equation} \label{jb}\begin{array}{ccc}ds_g^2 & =&  2\left(  Y(y) +{x}{} \right)dxdy \\ 
   ds_{\bar g}   & =&  -\frac{2(Y(y)+x)}{y^3}dxdy  + \frac{(Y(y)+x)^2}{y^4}dy^2.
 \end{array}\end{equation}

Arguing as above, we obtain that the condition \eqref{en1t}  is equivalent to a certain system of  6 PDEs  in the unknown functions $v_1, v_2, Y$.  
Solving the first 5 PDEs with respect to the derivatives of the unknown functions and
 substituting the solution in the remaining equation, we obtain that the system is equivalent to 
   \begin{equation} \label{jjjj} \left. \begin{array}{ccc} 
      {\frac {\partial v_1}{\partial x}}  & = & \frac{1}{2}y-\frac{3}{2}\\\frac{\partial v_1}{\partial y}        &=&\frac{1}{2}Y  +\frac{1}{2}x,
      \\v_2&= & y^2\\ Y'&  = & -\frac{1}{2}\frac{-5yY+3Y -5 yx+3x+2 v'_2 Y +2 v'_2   x+2v_1 }{v_2 } 
      \end{array}\right\}\end{equation}

    We see that the first three equations of \eqref{jjjj} are equivalent to $  v_2=  y^2 $,  $v_1 = 
    \left( \frac{1}{2}y
-\frac{3}{2} \right) x+\frac{1}{2} Y_1 $, where $Y_1'= Y$.  Substituting these in the last equation of \eqref{jjjj}, we obtain  the following linear ODE on $Y_1$: 
$$ y^2Y''_1 =\tfrac{1}{2}\left(y Y'_1  -3 Y'_1  
-{Y_1} \right).  
$$
The equation can be solved, the general solution is $$
  {Y_1}  = \left( y-3 \right) 
 \left(C_1+C_2\,\int_{y_0}^y {e^{\frac{3}{2\xi}}}\frac{\sqrt{|\xi|}}{\left( \xi-3 \right)^2}d\xi \right).  $$
Then, the function $Y= Y'_1$ is 
$$
{Y}  =   
 \left(C_2  {e^{\frac{3}{2y}}}\frac{\sqrt{|y|}}{\left( y-3 \right)}\right)+  
 \left(C_1+C_2\,\int_{y_0}^y {e^{\frac{3}{2\xi}}}\frac{\sqrt{|\xi|}}{\left( \xi-3 \right)^2}d\xi \right). 
$$
The assumption that the metric admits no  Killing vector field implies $C_2\ne 0$.  
 In view of  the coordinate change  $x_{\textrm{new}}  = x_{\textrm{old}}+ C_1$,     we can assume  $C_1=0$.

Then, the  components $v_1, v_2$ of the projective vector field are 
$$  v_2=  y^2 \ , \  \  v_1 = 
    \left( \frac{1}{2}y
-\frac{3}{2} \right) x+ C_2\,\frac{ y-3  }{2}
\int_{y_0}^y {e^{\frac{3}{2\xi}}}\frac{\sqrt{|\xi|}}{\left( \xi-3 \right)^2}d\xi  $$

 We see that the metric (after the appropriate coordinate change and  scaling)  and the projective vector field  $v$ 
  are as in case \ref{3a} of Theorem~\ref{main}.

    \subsection{ The matrix   of 
  $L_{v}$ is  as (b) in  (16)
   } \label{pb} 

Assume that in the basis $\{a  , \bar a\}$ the matrix of $L_v:\hat{\mathcal{A}}\to \hat{\mathcal{A}}$ is given by $\begin{pmatrix} \lambda  & -1 \\ 1 & \lambda  \end{pmatrix}.$

Without loss of generality, in view of Lemma~\ref{nond}, we can assume that $
a = \det(g)^{-2/3} \cdot  g\ , \  \   \ \bar a= \det(\bar g)^{-2/3} \cdot  \bar g$ for certain metrics $g, \bar g$ from the projective class of \eqref{eq1}. 
Then, as we explained in \S\ref{trick}, the condition 
$$L_v\begin{pmatrix} a \\ \bar a\end{pmatrix} =  \begin{pmatrix} \lambda  & -1 \\ 1 & \lambda    \end{pmatrix} \begin{pmatrix} a \\ \bar a\end{pmatrix}  $$
is equivalent to  the condition \eqref{cond3}.     

  As we explained in Appendix,
   in a neighborhood of almost every point the metrics $g$ and $\bar g$ have one of three normal forms. We will carefully consider all three cases. 

\subsubsection{ The  Liouville case }   \label{Lb}
   Assume the metrics $g$ and $\bar g$  have the Liouville form \eqref{liu}. 
    Then, the condition \eqref{cond3}  is  equivalent to a system of  6 PDE's 
   on the unknown functions $v_1(x, y)$, $v_2(x, y)$,  $X(x)$, and $Y(y)$.

  Solving the equations with respect to the derivatives, we obtain that the equations are equivalent to the following system of 6 PDEs in Frobenius form:  

  \begin{equation} \left.\label{e3}  
  \begin{array}{ccc}  {\frac {\partial v_1 }{\partial x{{}}}}{ } & =&
    X/2  -3/2\,\lambda    \\ {\frac {
\partial v_2  }{\partial y}} &=& 
Y/2  -3/2\,\lambda   \\ {\frac {
\partial v_1 }{\partial y{{ }}}} & =& 0 \\ {\frac {
\partial v_2  }{\partial x{{ }}}} & =& 0 \\ 
 X'{v_1} & =&{1+
 X^{2}}{  } \\ Y'v_2 & =& {
 {1+ Y^{2}}{{ }
 }} 
  \end{array} \right\}\end{equation} 

We see that the functions $v_1$ and $v_2$ are functions of one variable only\footnote{which was clear in  advance since the family $\hat{\mathcal{A}}$ determines the lines of the coordinates, see \cite{obata,CMH}; hence the coordinate lines must be preserved by the flow of $v$}, so that all the equations (\ref{e3}) are actually ODEs.  Moreover, the assumption that there exists no Killing vector field implies that $X$ and $Y$ are not constant. Then, in view of  the  first two equations of (\ref{e3}), the components 
 $v_1$ and $v_2$ are not zero almost everywhere.  Then, without loss of generality we can assume $v_1\ne 0$, $v_2\ne 0$. Take   the new  coordinate system   $\bigl(x(x_{old}),  y_{new}(y_{old})\bigr)$  
given by   \begin{equation} dx =\frac{1}{v_1} dx_{old} \   \   dy =\frac{1}{v_2} dy_{old}.\label{dx} \end{equation} In these new coordinates the last two equations of (\ref{e3}) are ${\dot  X}=   1 + X^2$,
 ${\dot  Y}=   1 + Y^2$ implying\begin{equation}X(x)  =  \tan(x + \const_1)\,,  \ \ Y(y)  =  \tan(y_{new} + \const_2)\label{dx1} .\end{equation}
  Since the      formulas (\ref{dx}) define $x$ and $y$ up to addition of a constant, we can (and will)  assume without loss of generality 
   that $\const_1=\const_2=0$. 
Now consider the first and the second   equations of  (\ref{e3}). In the new coordinates they  are  $\dot v_1   =
    \bigl( \tan(x)/2  -3/2\,\lambda\bigr)v_1,\ $
$\dot v_2   =
    \bigl( \tan(x)/2  -3/2\,\lambda\bigr)v_2.$
     Solving them  we obtain \begin{equation}\label{dx2b}  {v_1}  ={\frac {{ C_1}  {e^{-3/2\,\lambda\,x{{}}}}}{\sqrt {\cos \left( x{{}}
 \right) }}} \   , \  \   {v_2}  ={\frac {{ C_2}  {e^{-3/2\,\lambda\,y{{}}}}}{\sqrt {\cos \left( y{{}}
 \right) }}}. 
\end{equation}   
Finally,  combining (\ref{dx} -- \ref{dx2b}) we obtain  (after the appropriate coordinate change and the scaling) that 
the metric $g$ has the form 
\begin{equation} \label{fuck} ds^2_g= (X-Y) \left(dx_{old}^2 \pm dy_{old}^2\right) = (\tan(x)- \tan(y)) 
\left(  {\frac {{ C_1^2}  {e^{-3\,\lambda\,x}}d x^2}{ {\cos \left( x
 \right) }}}\pm {\frac {{ C_2^2}  {e^{-3\,\lambda\,y}}d y^2}{ {\cos \left( y
 \right) }}}\right), \end{equation}
and the projective vector field $v $ is $\frac{\partial }{\partial x } + \frac{\partial }{\partial y }$.
It is easy to see that if $\lambda=0$ and $C_1=\pm C_2$, then the metric admits a Killing vector field. Indeed,  
because of scaling, it is sufficient to show this for $C_1=C_2=1$. If the sign  ``$\pm$" in  \eqref{fuck} is ``$-$", then  the vector field $\frac{\cos(x)}{\sin\left(\tfrac{1}{2}(x-y) \right) }\frac{\partial}{\partial x}  -\frac{\cos(y)}{\sin\left(\tfrac{1}{2}(x-y) \right) }\frac{\partial}{\partial y}$ is a Killing
one  for  $g$. If the sign  ``$\pm$" in  \eqref{fuck} is ``$+$", then  
 the vector field $\frac{\cos(x)}{\cos\left(\tfrac{1}{2}(x-y) \right) }\frac{\partial}{\partial x}  +\frac{\cos(y)}{\cos\left(\tfrac{1}{2}(x-y) \right) }\frac{\partial}{\partial y}$ is a Killing
 one  for  $g$.

 We see that the metric   and the projective vector field  $v$ 
  are as in case \ref{1b} of Theorem~\ref{main}.

  \subsubsection{The Complex-Liouville case } \label{CLb} 

   Assume the metric $g$ and $\bar g$  have the Complex-Liouville form (\ref{gcomplex}). 
   Arguing as above, we obtain that   the conditions  
  \eqref{cond3} are equivalent to a certain 
  system  of   6 PDE of the first order. 
  Solving this system with respect to the first derivatives, and using  the Cauchy-Riemann conditions for the holomorphic function $h$, we obtain  that the system is equivalent to the  system 
  \begin{equation} \label{pp} 
  \left. \begin{array}{rcc} 
   v^2_y = {v^1_x}  &=&\frac{1}{2}{ \Re(h)} -\frac{3}{2}\lambda\\ 
-{v^1_y}  =v^2_x  &=&\frac{1}{2}{ \Im(h)}\\ 
 { \Re}(h)_x= 
 { \Im}(h)_y
  & =& {\frac {-{ v^1}   { \Im}^2(h)+2{ \Re(h)} { \Im(h)} {v^2} +{v^1} +{ v^1}   { \Re^2 (h)} 
}{ \left( { v^2}   \right) ^{2}+
 \left(v^1\right)^{2}}}  \\
  -\Re(h)_y= { \Im(h)}_x & =& {\frac {{ v^2}  { \Im^2(h)}
 -{v^2}-{v^2}  { \Re^2(h)}+2{v^1}{ \Re}(h) { \Im(h)}
 }{ \left( {v^2}\right)^{2}+ \left( {v^1}   \right) ^{2}}}
\end{array} \right\}\end{equation}

From the first two equations of \eqref{pp}  we see that the function $V:= v_1 + i\cdot v_2$ is a holomorphic function of the variable $z:= x +i\cdot y$. It is easy to check that the last  two equations of \eqref{pp} are equivalent 
to \begin{equation}  Vh_z= h^2 +1 \label{ppp}, \end{equation} and the first two equations of   \eqref{pp} are equivalent to \begin{equation}V_z= \frac{1}{2} h - \frac{3}{2}\lambda\label{pppp} \end{equation}
 (where $h_z$, $V_z$ are  the derivatives  of $h, V$ with respect to $z$). Thus, the equations \eqref{pp} are direct analogues of  (\ref{e3}). 
After the coordinate change   
\begin{equation} dz_{\textrm{new}} =\frac{1}{V} dz_{\textrm{old}} \label{dz},  \end{equation}   
the equation  (\ref{ppp}) reads $\frac{dh}{dz_{\textrm{new}}}=   1 + h^2$  implying\begin{equation}h(z_{\textrm{new}})  =  \tan(z_{\textrm{new}} + \const)\, \label{dz1} .\end{equation}
  Since the     formulas (\ref{dz}) define $z_{\textrm{new}}$ 
  up to addition of a constant, we can (and will)  assume without loss of generality 
   that $\const=0$. In this new coordinate the vector field $v$ is $V\frac{\partial }{\partial z} + \bar V \frac{\partial }{\partial \bar z}$. 

Now consider the equation (\ref{pppp}). In the new coordinates it reads \begin{equation}\label{dz2}{\frac {d V }{d z_{\textrm{new}}}}{ }  =
    \bigl( \tan(z_{\textrm{new}})/2  -\frac{3}{2}\,\lambda\bigr)V \  \  \  \textrm{implying} \ \  \ 
      {V}  ={\frac {{ C}  {e^{-3/2\,\lambda\,z_{\textrm{new}}}}}{\sqrt {\cos \left( z_{{\textrm{new}}}
 \right) }}}. 
\end{equation}  
Finally,  combining (\ref{dz} -- \ref{dz2}) we obtain that (after the appropriate scaling)
the metrics $g$  and $\bar g$ have  the form 
\begin{equation*}  
\begin{array}{ccl} ds^2_g & = &     \frac{1}{4} \left(\tan(z_{\textrm{new}})- \tan(\bar z_{\textrm{new}})\right) 
\left(  {\frac {{ \bar C}  {e^{-3\,\lambda\,\bar z_{\textrm{new}}}}d \bar z_{\textrm{new}}^2}{ {\cos \left( \bar z_{\textrm{new}}
 \right) }}}-  {\frac {{ C}  {e^{-3\,\lambda\,z_{\textrm{new}}}}d z_{\textrm{new}}^2}{ {\cos \left( z_{\textrm{new}}
 \right) }}}\right),\\ 
  ds^2_{\bar g} & = &  \frac{1}{4} \left(\cotan(\bar z_{\textrm{new}})- \cotan(z_{\textrm{new}})\right) 
\left(  {\frac {{ \bar C}  {e^{-3\,\lambda\,\bar z_{\textrm{new}}}}d \bar z_{\textrm{new}}^2}{ {\sin \left( \bar z_{\textrm{new}}
 \right) }}}-  {\frac {{ C}  {e^{-3\,\lambda\,z_{\textrm{new}}}}d z_{\textrm{new}}^2}{ {\sin \left( z_{\textrm{new}}
 \right) }}}\right), \end{array}\end{equation*} 
and the projective vector field $v $ is $\frac{\partial }{\partial z_{\textrm{new}} } + \frac{\partial }{\partial \bar z_{\textrm{new}} }= \frac{\partial }{\partial x_\textrm{new}}$.  
  It is easy to see that if «$\lambda=0$ and $C\in \mathbb{R}$,  then the metric admits a Killing vector field.  Indeed, it is sufficient to consider $C=1$. For this case,  the following 
   vector field is a Killing one:
  $
  \sin \left( x\right) \frac{\partial }{\partial x} +\frac {\cos \left( x
 \right) \cosh \left( y \right) }{\sinh \left( y \right) }\frac{\partial }{\partial y}.
$

 We see that the metric   and the projective vector field  $v$ 
  are as in case \ref{2b} of Theorem~\ref{main}.

 \subsubsection{The  Jordan-block case} \label{JBb_new} 

  Assume the metrics $g$ and $\bar g$ are given by the matrices \eqref{jb}.
Arguing as above, we obtain that the condition \eqref{cond3}  is equivalent to a certain system of  6 PDEs in the unknown functions $v_1, v_2, Y$.  
Solving the first 5 PDEs with respect to the derivatives of the unknown function,  and 
substituting the solution in the remaining equation, we obtain that the system is equivalent to 
   \begin{equation} \label{jj} \left. \begin{array}{ccc} 
      {\frac {\partial v_1}{\partial x}}  & = & \frac{1}{2}y-\frac{3}{2}\lambda\\\frac{\partial v_1}{\partial y}        &=&\frac{1}{2}Y  +\frac{1}{2}x,
      \\v_2&= & y^2+1\\ Y'&  = & -\frac{1}{2}\frac{-5yY+3\lambda Y -5 yx+3 \lambda x+2 v'_2 Y +2 v'_2   x+2v_1 }{v_2 } 
\end{array}\right\}\end{equation}

    We see that the first three equations of \eqref{jj} are equivalent to $$v_2=  y^2+1 \ ,  \  v_1 = 
    \left( \frac{1}{2}y
-\frac{3}{2} \lambda \right) x+\frac{1}{2} Y_1(y), \textrm{  where $Y'_1=  Y(y)  $ } $$  Substituting these in the last equation of \eqref{jj}, we obtain  the following linear ODE on $Y_1$: 
$$ (1+y^2)Y''_1 =\frac{1}{2}\left(y Y'_1  -3 \lambda Y'_1  
-{Y_1} \right).  
$$
The equation can be solved, the general solution is $$
  {Y_1}  = \left( y-3\lambda \right) 
 \left(C_1+C_2\,\int_{y_0}^y {e^{-\frac{3}{2}\lambda \arctan(\xi)}}\frac{\sqrt[4]{\xi^2+1}}{\left( \xi-3\lambda \right)^2}d\xi \right).  $$
Then, the function $Y= Y'_1$ is 
$$
{Y}  =    
 C_1+C_2\,\int_{y_0}^y {e^{-\frac{3}{2}\lambda \arctan(\xi)}}\frac{\sqrt[4]{\xi^2+1}}{\left( \xi-3\lambda \right)^2}d\xi +\left( y-3\lambda \right) 
 \left(C_2\,{e^{-\frac{3}{2}\lambda \arctan(y)}}\frac{\sqrt[4]{y^2+1}}{\left( y-3\lambda \right)^2} \right) 
$$
and the  components $v_1, v_2$ of the projective vector field are 
$$  v_2=  y^2 +1 \ , \  \  v_1 =  
    \left( \frac{1}{2}y
-\frac{3}{2}\lambda  \right) x+ \frac{ y-3\lambda  }{2}
\left(C_1+C_2\,\int_{y_0}^y {e^{-\frac{3}{2}\lambda \arctan(\xi)}}\frac{\sqrt[4]{\xi^2+1}}{\left( \xi-3\lambda \right)^2}d\xi \right). $$

 We see that, after an appropriate coordinate change and    scaling,   
  the metric   and the projective vector field  $v$ 
  are as in case \ref{3b} of Theorem~\ref{main}.

  \subsection{ The matrix  of $L_{v}$ is as (c) in (16) }  \label{pc} 
Assume that in the basis $\{a  , \bar a\}$ the matrix of $L_v:\hat{\mathcal{A}}\to \hat{\mathcal{A}}$ is given by $\begin{pmatrix} \lambda  &  \\  & 1  \end{pmatrix},$ where $\lambda \in (-\infty, -1] \cup [1, +\infty)$.

Without loss of generality, in view of Lemma~\ref{nond}, we can assume that $
a = \det(g)^{-2/3} \cdot  g\ , \  
\   \ \bar a= \det(\bar g)^{-2/3} \cdot  \bar g$ for certain metrics $g, \bar g$ from the projective class of \eqref{eq1}. 
Then, as we explained in  \S\ref{trick}, the condition 
$$L_v\begin{pmatrix}  a  \\ \bar a\end{pmatrix} =  \begin{pmatrix} \lambda  &  \\  & 1    \end{pmatrix} \begin{pmatrix} a \\ \bar a\end{pmatrix}  $$
is equivalent to  the condition 
$
  L_{v} a =  \lambda  a 
   $, $  L_{v} \bar a =  \bar a
 ,$ which is  equivalent to the condition  \[
L_vg-\frac{2}{3}\text{trace}_g(L_vg)g-\lambda g=0\,,\quad
L_v\bar{g}-\frac{2}{3}\text{trace}_{\bar{g}}(L_v\bar{g})\bar g -\bar{g}=0
,\] which is equivalent to  the condition
     \begin{equation}\label{homot} 
  L_{v} g = 
  -\frac{\lambda}{3}g \, , \ \ \  L_{v} \bar g = -\frac{1}{3} \bar g.
  \end{equation}

  As  explained in  the appendix,
   in a neighborhood of almost every point the metrics $g$ and $\bar g$ have one of three normal forms. We will carefully consider all three cases. 

   \begin{Rem} \label{addition}
   We will also see that $\lambda \ne 1$. This will imply that  if  two nonproportional projectively equivalent metrics $g$ and $\bar g$ have $L_v g =\lambda \cdot  g $ and $ L_v \bar g=\lambda\cdot  \bar g$
   for a certain $v \not\equiv  0$, then the metrics admit a Killing vector field, which will be used in the proof of Theorem~\ref{main2}. 
   \end{Rem}

  \subsubsection{The Liouville case} \label{Lc}

  We assume that the metrics $g$ and $\bar g$ are given by \eqref{liu}. Then, the condition \eqref{homot} is equivalent to a system of 6  PDEs in $v_1,  v_2, X, Y$. 

  Solving these equations  with respect to derivatives,  we obtain 

  \begin{equation} \label{long} 
\begin{array}{lcr}
       \frac{\partial v_1 }{\partial x}   =  -\lambda-1/2,  &  \frac{\partial v_1}{\partial y}  =  0,  &   X'v_1   =  -X  \left( -1+\lambda
 \right) \\ \frac {\partial v_2 }{\partial y} =  -\lambda-1/2,     &  \frac {\partial v_2 }{\partial x}=0, & 
Y' v_2  =  -Y   \left( -1+\lambda \right)  \end{array} . \end{equation}

  This system of equations can be easily solved (we recall that by assumption $|\lambda|\ge 1$).  
  If $\lambda = 1$, at least one of the functions $X, Y$ is a constant implying the existence of a Killing vector field as  promised in Remark \ref{addition}.   
  For other $\lambda$,    the solution is $\bigl($up to the coordinate change $(x_{new}, y_{new}) =(x_{old} + \const_1, y_{old}+ \const_2)$ $\bigr)$ 
  \begin{equation*} 
\begin{array}{c}v_1= -\left(\frac{1}{2}+ \lambda\right) x , \  \ X= C_1\, x^{2\tfrac{ \lambda-1}{1+ 2\lambda}} \\  
v_2= -\left(\frac{1}{2}+ \lambda\right) y , \  \ Y= C_2\, y^{2\tfrac{\lambda-1}{1+ 2\lambda}}, \end{array} \end{equation*}
and the corresponding   $g$  and $v$, after dividing $v$  by $-\left(\tfrac{1}{2}+ \lambda\right)$,  are 
 $$\left(C_1\, x^{2\tfrac{\lambda-1}{1+ 2\lambda}}  - C_2\, y^{2\tfrac{\lambda-1}{1+ 2\lambda}}\right) (dx^2 + \varepsilon dy^2)\ , \ \
   x \frac{\partial }{\partial x} +  y \frac{\partial }{\partial y}. $$

 We see that after the coordinate change $(x_{\textrm{old}}=e^{x}, y_{\textrm{old}}=e^{y})$,  after an appropriate scaling, and after denoting $\tfrac{2(\lambda -1)}{2\lambda +1 }$  by $\nu$,   
  the metric   and the projective vector field  $v$ 
  are as in case \ref{1c}   of Theorem~\ref{main}. Note that in 
  the case  $\nu =2$, 
  $ C_1=- \varepsilon  C_2$ the  metric $g$ has a  constant curvature (and, therefore, a Killing vector field). 
  Since $\lambda\in (-\infty, -1]\cup (1,+\infty)$, we have  $\nu \in (0, 4],$  $ \nu\ne 1$.   
    Since $\lambda \ne 1$, then $\nu\ne 0$.

  \subsubsection{The Complex-Liouville case} \label{CLc}

   Assume that 
   $g$,  $\bar g$ are  as in (\ref{gcomplex}).
Arguing as above, we obtain that the equations  \eqref{homot} are equivalent to a system of 6 PDEs which can be written  as
\begin{eqnarray*}
 v_1+ i\cdot  v_2&=& -(\lambda+ 1/2)\cdot z +  \const   \\
     \frac{\partial h}{\partial z}&=& \frac{h}{z} \cdot \frac{2(\lambda-1)}{1+ 2\lambda}.\end{eqnarray*} 
 The system can be easily  solved. If $\lambda = 1$,  the function $h$ is a constant implying the existence of a Killing vector field as we promised in Remark \ref{addition}.

  If $\lambda \ne 1$, then, in view of the coordinate change  $x_{new} =  x_{old} + \const_1,\ y_{new} =  y_{old} + \const_1$ we can take   $\const=0$. Then, the solution is $h= C\cdot z^{2\frac{(\lambda-1)}{1+ 2\lambda}}$. 
 Then, the metrics $g$ and $\bar g$ are as in \eqref{complexification} with this  function $h$, 
 and  the projective vector field  $v$ is $ x\frac{\partial }{ \partial x }$. 

 We see that after the coordinate change $z_{\textrm{old}}=e^{z}$,
 after an appropriate scaling, and  after denoting $\tfrac{2(\lambda -1)}{2\lambda +1 }$  by $\nu$,   
  the metric   and the projective vector field  $v$ 
  are as in case \ref{2c} of Theorem~\ref{main}.  Since $\lambda\in (-\infty, -1]\cup (1,+\infty)$, we have  $\nu \in (0, 4]$, $\nu\ne 1$,   as we assumed in Theorem~\ref{main}.

 \subsubsection{The Jordan-block  case} \label{JBc} 
 Assume the metrics $g$ and $\bar g$ are given  by \eqref{jb}.  
Arguing as above, we obtain that the condition \eqref{homot}  is equivalent to a certain system of  6 PDE on the unknown functions $v_1, v_2, Y$. 
Solving this system with respect to the derivatives of the functions we see that $\frac{\partial v_2}{\partial x} =0$ , $\frac{\partial v_1}{\partial y} =0$
implying that $v_1$  is a function of $x$ and  $v_2$ is a function of $y$ only, and that 
    the system is equivalent to  

    \begin{equation} \label{60} 
    \left.\begin{array}{ccc}
     v'_1&=&-\lambda-\frac{1}{2}\\v_2 &=&- \left( 
\lambda-1 \right) y \\
     Y'{v_2 } &=&-\frac{1}{2}\,\left(2 v'_2 Y  +2\,   v'_2 x+4\lambda Y -Y +4\lambda x-x+2v_1 \right)
  \end{array}\right\} \end{equation} 

From the first equation of \eqref{60}  we see that $v_1=-\lambda x -\frac{x}{2} +C$. Without loss of generality we can take $C=0$. Substituting the expressions for $v_1, v_2$ in the last equation of \eqref{60}, we obtain (we can assume $y>0$ since it can be achieved by a coordinate change)
$$Y'  =\frac {(2\lambda+1)\,Y }{\left( 2\lambda-2 \right) y}. 
$$
 Solving this equation, we obtain 
 \begin{equation}  Y  =y^{{\frac{2\,\lambda+1}{2\lambda-2}}}{C_1} . \end{equation}

 We see that after an appropriate scaling and 
  after denoting $\tfrac{2(\lambda -1)}{2\lambda +1 }$  by $\eta$,
  the metric   and the projective vector field  $v$ 
  are as in case \ref{3c} (for $\eta \ne \tfrac{1}{2}$) or as in case \ref{3d} (for $\eta  =\tfrac{1}{2}$) of Theorem~\ref{main}.

\section{ Proof of Theorem 2}

\subsection{ In the cases 1a -- 3c,   it is sufficient to prove that $\mathcal{A}$ is precisely 
two-dimensional. In the case 3d,  it is sufficient to prove that  $\mathcal{A}$ is precisely 
three-dimensional.} \label{4.1}

Within this paragraph  we assume that the metric $g$ is one of the metrics from  Theorem \ref{main}. We additionally assume that it   admits no Killing
vector field. Let us explain why, in order   to prove Theorem~\ref{main2},  it is sufficient 
 to show that the space $\mathcal{A}$ of solutions of \eqref{lin1} is as in the title of this paragraph.

  For every metric   $g$  from Theorem 
  \ref{main}, consider its canonically projectively equivalent metric $\bar g$ given by the appropriate formula from (\ref{E1} -- \ref{E3}).  By definition, the metrics $g$ and $\bar g$ have the same  projective connection. Then, 
   $a = g/\det(g)^{2/3}$ and   $\bar a = \bar g/\det(\bar g)^{2/3}$ lie in the space $\mathcal{A}$ corresponding to the metric $g$. 

    Therefore, every  linear combination $\alpha \cdot a + \beta \cdot \bar a $ is also an element of $\mathcal{A}$. Comparing the definition of $G(g, \bar g)$ with the formulas in Remark \ref{rem3}, we see that  $G(g, \bar g)$ is precisely the set of  metrics corresponding to the solutions of the form $\alpha \cdot a + \beta \cdot \bar a $.    In particular, all metrics from $G(g, \bar g) $  lie in the projective class of $g$. 

    Thus, in order to show that the projective class of the metrics $g$ from cases \ref{1a} -- \ref{3c} of Theorem~\ref{main} coincides with  $G(g, \bar g) $, it is sufficient to show that  $\mathcal{A}$ coincides with the set of linear combinations of $a$ and $\bar a$,  i.e.,  is  two-dimensional.

 Now, let us consider 
  the metric \ref{3d}. In this case, the space  $\mathcal{A}$    is at least  three-dimensional. Indeed, the solutions    $a = g/\det(g)^{2/3}$,  
  $\bar a = \bar g/\det(\bar g)^{2/3}$, and $\tilde a = \tilde g/\det(\tilde g)^{2/3}$ are linearly independent.
    Clearly, the metrics 
  corresponding to the linear combinations of these solutions are precisely the metrics from $G[g,\bar g, \tilde g]$. Hence, if  the space of $\mathcal{A}$ is precisely  three-dimensional,    the projective class coincides with $G[g,\bar g, \tilde g]$.

\subsection{ Schema of the  proof } 

Theorem 1 gives us 10 explicit formulas for the  metric $g$ 
and,  therefore,   for the coefficients $K_i$ of the equation \eqref{lin1}.     Our goal is to show that in the first 9 cases the space $\mathcal{A}$ 
of the solutions  of  \eqref{lin1} is at most  two-dimensional, and in the last case the space $\mathcal{A}$  is  at most three-dimensional.

There exists a highly computational method to do it: indeed, the system \eqref{lin1} is linear
 and of finite type. Then, the standard prolongation-projection method gives us an algorithm which calculates the dimension of the space of solutions.

Unfortunately, this method is too hard from the viewpoint of calculations,  at least if one does the calculations straightforwardly: indeed, in order to implement the algorithm, 
  one needs to differentiate the entries of the metric  $7$ times, and then calculate the rank of an  
 $18\times 16$ matrix.

 It is  still possible to  do it  with the help of computer algebra packages.
  Recently Kruglikov \cite{kruglikov} and, independently,  
  Bryant,  Eastwood and Dunajski  \cite{duna} used Mathematica$^\registered$ and Maple$^\registered$
  (and also quite advanced theory) to construct   curvature   invariants  such  that  if they do not vanish  
  the dimension of $\mathcal{A}$ is at most $2$.   But their invariants   are still too complicated, and there is no hope    to calculate them for our metrics  without  using a computer (though one can easily do it by  computer).

 In order to give  a proof which is much easier from a computational point of view, and which could be done by a human,  we use the existence of the projective vector field to reduce the problem to more simple systems of PDEs. We consider {\bf three cases}. 

 The {\bf first case} corresponds to the metrics \ref{3a}, \ref{3b}, \ref{3c},  \ref{3d}.  
 In these cases,  the general form of the metric is very simple and 
 one can actually  do the prolongation-projection algorithm by hand and without using the existence of the projective vector field, see \S\ref{case1}. After a few  steps (we actually use short-cuts in the paper),    we obtain  the dimension of $\mathcal{A}$. 

The {\bf second  and the third cases}  corresponds to all  other metrics.  We assume that  $\dim(\mathcal{A})=3$ and find a contradiction. (The case $\dim(\mathcal{A})\ge 4$
 is not possible because by Theorem \ref{main3} the metrics \ref{1a} -- \ref{2c} admit no Killing vector field. We will not use Theorem 2 in the proof that the metrics \ref{1a} -- \ref{2c} admit no Killing vector, so no logical loop  appears). 
In order to do it, let us take a  basis   $\{a,\bar a,  \hat a \}$ such that the matrix of 
 $ L_v$ is one of \eqref{mat3bis}, where $\textrm{ \Large $A$}$ is a $(2\times 2)$- matrix given by \eqref{mat}
 \begin{equation}\label{mat3bis} 
   \begin{pmatrix} \mu &  \\ &  \textrm{ \Large $A$}  \end{pmatrix}, \ 
    \ \begin{pmatrix} 1 & 1& \\ &  1 &1 \\ &&1 \end{pmatrix}\, .  
   \end{equation}

  The {\bf second case} corresponds to the first  matrix of 
  \eqref{mat3bis}.   We will  find a contradiction using the following   trivial observation from linear algebra: if for a 
  $(3\times 3)$ matrix $m= (m_{ij}) $ with $det(m)\ne 0$ 
  \begin{equation} \label{m}\left. 
\begin{array}{ccc}   m_{11} a_{11} + m_{12} a_{12} + m_{13} a_{22}  &= & 0 \\ 
m_{21} a_{11} + m_{22} a_{12} + m_{23} a_{22}  &= & 0 \\
 m_{31} a_{11} + m_{32} a_{12} + m_{33} a_{22}  &= & 0,   \end{array} \right\}\end{equation} 
  then $a_{ij}=0$.

  Let us explain how the assumptions of the {\bf second case} allow us  to construct such equations in $a_{ij}$. 
  Since the matrix of $L_v$ is the first matrix of \eqref{mat3bis},   
   we have \begin{equation} L_v\begin{pmatrix} a \\\bar a \\ \hat a\end{pmatrix} = \begin{pmatrix} \mu &  \\ &  \textrm{ \Large $A$} \end{pmatrix}\begin{pmatrix} a \\\bar a \\ \hat a\end{pmatrix}.  \label{tmp} \end{equation} The last two equations of \eqref{tmp}  are equations $L_v \begin{pmatrix} \bar a \\ \hat a\end{pmatrix} =     \textrm{ \Large $A$}  \begin{pmatrix} \bar a \\ \hat a\end{pmatrix}$.  We solved them in the proof of Theorem \ref{main}. In a certain coordinate system (in a neighborhood of almost every point) the metric $g= \bar a/det(\bar a)^2$ is as in Theorem \ref{1} after a possible scaling.  We have 6 (explicit) possibilities \ref{1a} -- \ref{2c} for the metric and therefore 6 (explicit)  possibilities for the coefficients of the equation \eqref{lin1}.

Let us now  pass to the coordinate system  such that the projective vector field is $\frac{\partial}{\partial x}$. In cases \ref{2a} -- \ref{2c}, 
we are already in such a coordinate system, in the cases  \ref{1a} -- \ref{1c} we use the coordinate change $x_{new}= \frac{x_{old} + y_{old}}{2}$, $y_{new}=\frac{ x_{old}- y_{old}}{2}$. 
In this coordinate system,   the coefficients 
 $K_0,...,K_3$ of the projective connection are independent of $x$ and direct calculations show that they are given  by simple formulas, see the beginning of \S\ref{case2}.   
   We see that  the first equation of \eqref{tmp} is  
   $\frac{\partial  a}{\partial x}  = \mu \cdot  a $ implying 
  \begin{equation}  a = e^{\mu x} \begin{pmatrix} a_{11}(y) & a_{12}(y) \\ a_{12}(y) & a_{22}(y)\end{pmatrix}.\label{hom} \end{equation}  Substituting \eqref{hom}  in 
  \eqref{lin1}, we obtain    one homogeneous linear   equation 
    and 3 linear ODEs in the three unknown functions $a_{ij}(y)$, see  \S~\ref{anal} for the precise formulas. This linear  equation (first equation  of \eqref{s1})   will play the role of  the first equation of \eqref{m}.  

      It is possible to explicitly solve the   above mentioned 
   3 ODEs  with respect to derivatives, see  \eqref{s1}. 
  Differentiating the first equation  of \eqref{s1} with respect to $y$, and substituting the derivatives of $a_{ij}$ from the other three equations of 
  \eqref{s1} inside, we obtain one more linear equation  on $a_{ij}$. This equation  will play the role of the second equation of \eqref{m}.  Repeating the procedure with this new equation, we obtain  the third linear equation  on $a_{ij}$, which will be  the third equation of \eqref{m}.  Direct calculations  show that the determinant of the correspondent $(3\times 3)$-matrix  $(m_{ij})$ is not zero  implying $a\equiv 0$.
  We obtain a contradiction with  the assumption that $\{a, \bar a , \hat a\}$ is a basis.

   This described procedure is not very complicated computationally  
   (all formulas that appear have less than 50 terms, i.e., one can do all calculations by hand,  and standard computer algebra  packages, say Maple$^\registered$ or Mathematica$^\registered$,  need less then 10  seconds for all the calculation.)

  Let us also note that
  the proof for  the cases  \ref{1a}, \ref{1b}, \ref{1c}  implies 
  the proof for  the cases  \ref{2a}, \ref{2b}, \ref{2c} (so we need to do the calculation for the three cases \ref{1a}, \ref{1b}, \ref{1c} only).
   Indeed, the formulas for (the components of)  the metrics  \ref{1a}, \ref{1b}, \ref{1c} are real-analytic, and we can allow $x$ to be a complex variable and $y$ to be its conjugate, since it  changes neither  differentiation nor algebraic operations with the (components of the) metrics.  After this change  the metrics \ref{1a}, \ref{1b}, \ref{1c}  become, up to a multiplication by a constant, the metrics     \ref{2a}, \ref{2b}, \ref{2c}, and therefore our proof that the metrics \ref{1a}, \ref{1b}, \ref{1c}  have two-dimensional $\mathcal{A}$, which uses only algebraic operations and differentiation,    is also a proof for the cases 
    \ref{2a}, \ref{2b}, \ref{2c}. 

   The {\bf third case} corresponds to the second  matrix of 
  \eqref{mat3bis}.    We will  find a contradiction 
    using  the following fact  from linear algebra:  if 
    \begin{equation}\label{tmp1} \left.
\begin{array}{ccc}   m_{11} a_{11} + m_{12} a_{12} + m_{13} a_{22}  &= & b_1 \\ 
m_{21} a_{11} + m_{22} a_{12} + m_{23} a_{22}  &= & b_2 \\
 m_{31} a_{11} + m_{32} a_{12} + m_{33} a_{22}  &= & b_3 ,   \end{array} \right\}\end{equation} then the following two statements  are contradictive: 

det$\left(
\begin{array}{ccc}    m_{11} &  m_{12}  &  m_{13}  \\ 
m_{21}   &  m_{22} &   m_{23}   \\
 m_{31} &  m_{32} &  m_{33}   \end{array} \right)=0 \ ,  $ 
 \ \ \  \ \  and \  \ \ \ \ det$\left(
\begin{array}{ccc}  b_1  &  m_{12} &    m_{13}  \\ 
b_{2} & m_{22}   &    m_{23}   \\
 b_{3} & m_{32} &   m_{33}   \end{array} \right)\ne 0.   $

  The way to construct  equations \eqref{tmp1} are similar to that we use in the {\bf second case}.
   Since the matrix of $L_v$ is the second matrix of \eqref{mat3bis}, the Lie derivatives of the basis elements $a, \bar a , \hat a$ are given by 3 matrix equations  
     \begin{equation} \label{maned} 
     L_v \begin{pmatrix} a\\ \bar a   \\ \hat a\end{pmatrix} = \begin{pmatrix} 1  &1 & \\ & 1& 1 \\ && 1\end{pmatrix} \begin{pmatrix} a\\ \bar a   \\ \hat a\end{pmatrix}= 
     \begin{pmatrix} a  &+& \bar a & & \\ & &  \bar a & + & \hat a  \\ &&&& \hat a \end{pmatrix}  .\end{equation}       
 We see that the last two  equations  of \eqref{maned} are the equations \eqref{17}. We solved them  in \S\ref{J}, see \S\S \ref{La}, \ref{CLa} there.   Then,  
    without loss of generality  we can assume that 
    the metric $g= \bar a/\det(\bar a)^2$   and the projective field $v$ are  
      as in cases  \ref{1a},  \ref{1b}  of Theorem \ref{main}.    We  again pass to the coordinates such that  the projective   vector field $v$  is $\frac{\partial }{\partial x}$: in the case \ref{1b}, we are already in these  coordinates, 
    in the case \ref{1a}, we will work in the coordinates  $x_{new}= \frac{x_{old} + y_{old}}{2}$, $y_{new}=\frac{ x_{old}- y_{old}}{2}$.  In these coordinates, the coefficients  $K_0, ..., K_3$  of the projective connection are  independent of   $y$, and the components of the Lie-derivative $L_va$  
    are       the $x-$derivatives of the components of $a$.  Then, the first equation of \eqref{maned} is 
  \begin{equation}  \frac{\partial  }{\partial x}  a= a + \bar a.\label{cb}\end{equation} This equation is actually a system of three equations, since $a$ is a symmetric   $2\times 2-$matrix. In this equation, 
   $\bar a$ is known: in view of Remark \ref{rem3}, it 
    is given  by $ g/\det(g)^{2/3}$,  and above we assumed that 
     $g$ is the metric     \ref{1a} from Theorem \ref{main}.
      Direct calculations shows that $a$ is given by  \eqref{bara}. Then,  \eqref{cb}   is a system of linear nonhomogeneous equations,  every  solution is  the  sum of a  partial solution  $P$ (for the metric \ref{1a}, a partial solution is \eqref{P}) and a solution of      the equation   $\frac{\partial a}{\partial x} =  a$, i.e., has the following form: \begin{equation}  \label{ansatz} 
e^x\cdot \begin{pmatrix} a_{11}(y)  &a_{12}(y) \\ a_{12}(y) & a_{22}(y)\end{pmatrix}  + P.\end{equation}

    Substituting the ansatz  \eqref{ansatz} in  the   equations \eqref{lin1}, we obtain one  nonhomogeneous linear equation (which will play the role of the first equation of \eqref{tmp1}),  
       and three nonhomogeneous linear  ODE of the first order 
     on the components $a_{ij}(y)$.  The ODE can be solved with respect to  the  derivatives of $a_{ij}$, see \eqref{s2}.

    Differentiating the above mentioned linear equation (which is the first equation of  \eqref{s2}) 
     with respect to $y$, and substituting the derivatives of $a_{ij}$ from the other equations of  \eqref{s2} inside, we obtain one more linear nonhomogeneous  equation  on $a_{ij}$. Repeating the procedure with the obtained equation, we obtain  the third linear   nonhomogeneous  equation  on $a_{ij}$. Thus, we have three nonhomogeneous  linear relations on three functions $a_{ij}$ as in \eqref{tmp1}.  If we   show that  these three nonhomogeneous  linear relations are not compatible, then the  dimension of $\mathcal{A}$ is at most 2. 

    Clearly,  the determinant of $(m_{ij})$ is zero,  since   $\hat a := \bar g/\det(\bar g)^{2/3} $ is a solution of the system \eqref{lin1} and of the equation  $ \frac{\partial}{\partial x} \hat a = \hat  a$,  and, therefore, gives us a solution of the homogeneous part of the above mentioned linear relations. 
     Direct calculations   show that   det$\left(
\begin{array}{ccc}    b_{1}  & m_{12} &   m_{13}  \\ b_{2} &
m_{22}   &     m_{23}   \\ b_{3} & 
 m_{32} &   m_{33}   \end{array} \right)\ne 0,   $  see \eqref{ror}. This gives us a contradiction which proves Theorem~\ref{main2}  for the metrics \ref{1a}, \ref{1b} from Theorem~\ref{main2}.

    Let us also note that, similar to the {\bf second case}, 
  the proof for  the case  \ref{1a}   implies 
  the proof for  the case  \ref{2a}. Indeed, the formulas for (the components of)  the metrics  \ref{1a}  are real-analytic, and we can allow $x$ to be a complex variable and $y$ to be its conjugate, since it  changes neither  differentiation nor algebraic operations with the (components of the) metric.  After this change  the metric \ref{1a}  becomes, up to a multiplication by a constant, the metric     \ref{2a},  and therefore our proof that the metric \ref{1a} has two-dimensional $\mathcal{A}$, which uses only algebraic operation and differentiation,    is also a proof for the case 
    \ref{2a}.

\subsection{ Calculations related to proof of Theorem 2  for  the metrics from cases
 1a, 1b, 1c from Theorem 1  assuming that the matrix of $L_v$ is as the first  matrix of (51) } \label{case2}  \label{anal} 

 For the metrics \ref{1a}, \ref{1b}, \ref{1c} from Theorem~\ref{main}, 
  the projective connections in the new coordinates  $x = \frac{x_{old} + y{old}}{2}$, 
  $y = \frac{x_{old} - y{old}}{2}$  are  respectively  given by 
 \begin{equation} \label{pcb}\begin{array}{ccccc}{y''}& = & \frac{ \left(   {e^{6\,y}}+{{c}}^{2}{e^{-6\,y}}+2\,{c}\,{} \right) }{8{c}{}y}&+&
 \frac {3 \left( 4\,y{c}\,{}
-{e^{6\,y}} +{{c}}^{2}{e^{-6\,y}}
 \right) {}}{8\,{c}\,{}\, y}y'
 \\  &+ & {\frac { \left( -2
\,{ c}\,{}+3\,{e^{6\,y}}+3\,{{c}}^{
2}e^{-6\,y} \right) { }}{8{ c}\,{}y}}(y')^2&-&
\,{\frac { \left( 12\,y{ c}\,{ }-{{c}}^{2}{e^{-6\,
y}}+{e^{6\,y}} \right) {}}{8{
c}\,{}y}}(y')^3\end{array} 
\end{equation} 

$$
\begin{array}{ccccc} 
y''&=&{\frac { {e^{6\,\lambda\,y}} + {c^2}{}{e^{-6\,\lambda\,y}}  +2\,c   \cos \left(2 y \right)  }{4c  \sin \left( 2 y \right) }}
&+&3\,{\frac {  -{e^{6\,\lambda\,y}}+
 {c}^{2}{e^{-6\,\lambda\,y}} +2c \,\lambda\, \sin \left( 2 y \right)  }{4 c\sin \left( 2 y \right) }} y'\\ &+&{\frac {  3\,{e^{6\,\lambda\,y}}+3\,{c}^{2}  {e^{-6\,\lambda\,y}}-2\,c \cos \left(2 y \right)  
}{4c\sin \left( 2 y \right) }} (y')^2&-&{\frac {{e^{6\,\lambda\,y}}-{c}^{2}+6\,\lambda\,c
 \sin \left( 2 y \right)  }{4 c \sin \left(2 y \right) }}(y')^3 \end{array}
$$

$$\begin{array}{ccc} 
y''&=& \frac{ -\lambda-\lambda\,{e^{4\,y}}+\lambda\,{e^{4\,y \left( \lambda-1 \right) }}{c}^{2}+\lambda\,{e^{4\,\lambda\,y}}{c}^{2}+\lambda\,{e^{2\,y \left( 2+\lambda \right) }}c-\lambda\,{e^{2\,y \left( -2+\lambda \right) }}c }{\left( 2 e^{2\lambda y}- 2\right)^2} \\ &+& \frac{ -4\,{e^{4\,\lambda\,y}}{c}^{2}+8\,c\,{e^{2\,\lambda\,y}}-4+3\,\lambda\,{e^{2\,y \left( 2+\lambda \right) }}c+\lambda\,{e^{4\,\lambda\,y}}{c}^{2}-2\,{e^{2\,\lambda\,y}}\lambda\,c+3\,\lambda\,{e^{2\,y \left( -2+\lambda \right) }}c-3\,\lambda\,{e^{4\,y \left( \lambda-1 \right) }}{c}^{2}+\lambda-3\,\lambda\,{e^{4\,y}} }{\left( 2 e^{2\lambda y}- 2\right)^2}y' \\ &+& \frac{ 3\,\lambda\,{e^{4\,y \left( \lambda-1 \right) }}{c}^{2}+3\,\lambda\,{e^{2\,y \left( 2+\lambda \right) }}c-\lambda\,{e^{4\,\lambda\,y}}{c}^{2}-3\,\lambda\,{e^{2\,y \left( -2+\lambda \right) }}c+\lambda-3\,\lambda\,{e^{4\,y}} }{\left( 2 e^{2\lambda y}- 2\right)^2} (y')^2\\ &+& \frac{ 4\,{e^{4\,\lambda\,y}}{c}^{2}-\lambda\,{e^{4\,y \left( \lambda-1 \right) }}{c}^{2}-8
\,c\,{e^{2\,\lambda\,y}}-\lambda-\lambda\,{e^{4\,y}}+4+\lambda\,{e^{2\,y \left( 2+\lambda \right) }}c-
\lambda\,{e^{4\,\lambda\,y}}{c}^{2}+2\,{e^{2\,\lambda\,y}}\lambda\,c+\lambda\,{e^{2\,y \left( -2+\lambda \right) }}c }{\left( 2 e^{2\lambda y}- 2\right)^2} (y')^3\end{array}$$

We see that all coefficients $K_i$ of the projective connection are independent of $x$ (which was clear in advance since the vector field $\frac{\partial }{\partial x}$ is projective).  
Substituting \eqref{hom} in  \eqref{lin1}, we obtain 
$$\left.\begin{array}{rcc} 
{ a_{11}}  \mu-\frac{2}{3}\,{ K_1}  { a_{11}}  +2\,{ K_0}  { a_{12}}&=&0 ,\\ { a'_{11}}  +2\,{ a_{12}}
  \mu-\frac{4}{3}\,{ K_2} 
   { 
a_{11}}  +\frac{2}{3}\,{ K_1} { a_{12}}  +2\,{ K_0} { a_{22}} &=&0,\\
2\,{ a'_{12}}  +{ a_{22}} \mu-2\,{ K_3}  { a_{11}}  -\frac{2}{3}\,{ K_2}  { a_{12}} +\frac{4}{3}\,{ K_1} 
{ a_{22}} &=&0,\\  
 { a'_{22}} -2\,{ K_3} { a_{12}} +\frac{2}{3}\,{ K_2}  
 { a_{22}} &=&0, 
\end{array}\right.$$
which is equivalent to the system 
 \begin{equation} 
\label{s1} \left. 
\begin{array}{ccc} 0&=&
{ a_{11}}  \mu-\frac{2}{3}\,{ K_1}  { a_{11}}  +2\,{ K_0}  { a_{12}} \\ { a'_{11}}  & = &  - 2\,{ a_{12}}
  \mu+\frac{4}{3}\,{ K_2} 
   { a_{11}}  -\frac{2}{3}\,{ K_1} { a_{12}}  -2\,{ K_0} { a_{22}} \\
{ a'_{12}} & =&  -\tfrac{1}{2}{ a_{22}} \mu+{ K_3}  { a_{11}}  +\frac{1}{3}\,{ K_2}  { a_{12}} -\frac{2}{3}\,{ K_1} 
{ a_{22}} \\  
 { a'_{22}}& =&  2\,{ K_3} { a_{12}} -\frac{2}{3}\,{ K_2}  
 { a_{22}}. 
\end{array}\right\}\end{equation}
The coefficients  $K_i$ of our connections are functions of  $y$ only.
 Differentiating the first equation of \eqref{s1} by $y$ and substituting the 
 values of $y$-derivative of $a_{ij}$ given by the last three equations, we 
  obtain the following equation as a differential consequence of the  equations \eqref{s1}: 
   if $a_{ij}$ satisfy \eqref{s1}, then they must satisfy the equation below. 
 \begin{equation*}
0= \frac{ {4}\mu{ K_2} -{2}{ K'_1}  +6{ K_0} { K_3} -{\frac{8}{3}}{ K_1} { K_2}   }
  {3} a_{11} + \frac{ {2}\mu K_1 +\frac{4}{3} K_1^{2}+6
{ K'_0}  -6{\mu}^{2}+2{ K_0}  { K_2}  
 }{3} { a_{12}}  -3\,\mu\,{ K_0}  { a_{22}}  
\end{equation*} 

Differentiating this equation 
  by $y$ and substituting the values of the $y$-derivative of $a_{ij}$ 
  from the last three equations of \eqref{s1}, we obtain another 
  linear homogeneous equation in 
  $a_{ij}$, whose coefficients are polynomial expressions in $K_i$ and their derivatives. 

Thus, we have three homogeneous linear equations  
in three unknown functions  $a_{ij}$, which must be satisfied if $a_{ij}$ satisfy \eqref{s1}. 
 The determinant of the corresponding $3\times 3-$ matrix is  given by \\
$
-2\,{\mu}^{6}+\frac{14}{3}\,{\mu}^{3} K_1    K_0
    K_2   +10\,  K'_0 {\mu}^{2} K_0    K_2   -{\frac {32}{9}}\,  K'_0  K_1^{2} K_0    K_2   -{\frac {40}{27}}\,\mu\, K_1^{3} K_0    K_2   +{\frac {64}{729}}\, K_1^{6}+6\,{\mu}^{2} K_0    K''_0   -10\,{\mu}^{4} K_0    K_2   -\frac{16}{3}\,{\mu}^{2}K_1^{2}K'_0-\frac{14}{3}\,{\mu}^{3} K_1   K'_0   -{\frac {64}{81}}\,K_1^{4} K_0    K_2   +{\frac {64}{27}}\, K_1^{3} K_0^{2} K_3   -{\frac {64}{81}}\, K_1^{3} K_0K'_1   +{\frac{16}{9}}\,K_0^{2}K_2^{2}K_1^{2}+\frac{8}{3}\,\mu\, K_1   (K'_0)^{2}+{\frac{40}{27}}\,\mu\, K_1^{3}K'_0   + \frac{8}{3}\,\mu\,K_1^{2}K_0^{2} K_3   -{\frac {32}{9}}\,\mu\,K_1^{2} K_0   K'_1   +{\frac {20}{3}}\,\mu\, K_0^{2} K_2  K'_1   -8\,\mu
\,  K_0^{3} K_2    K_3   +\frac{16}{3}\,{\mu}^{2}  K_1^{2} K_0
    K_2   +4\,{\mu}^{2}{
 K_1}    K_0   K'_1   -{\frac {32}{3}}\, K'_1  K_0^{3} K_3   -8\,{\mu}^{2}K_2^{2}K_0^{2}-12
\, K_0^{3}\mu\, K'_3   -4\, K_1   \mu\, K_0    K''_0   +{\frac {32}{3}}\, K_1 K'_0  K_0^{2} K_3   -{\frac {32}{9}}\, K_1  K'_0  K_0   K'_1   -\frac{4}{3}\,\mu\,  K'_0  K_1    K_0    K_2   -{\frac {32}{3}}\, K_1  K_0^{3} K_2    K_3   +4\, K_1   \mu\,  K_0^{2} K'_2   +4\, K_0^{2}\mu\,K''_1   +\frac{8}{3}\,\mu\, K_0^{2} K_2 ^{2} K_1   -16\,\mu\,  K'_0       K_0^{2} K_3   -\frac{8}{3}\,\mu\, K'_0   K_0    K'_1   +{\frac {32}{9}}\, K_1    K_0 ^{2} K_2    K'_1   +{\frac {16}{81}}\,\mu\, K_1 ^{5}+{\frac {64}{81}}\, K_1^{4} K'_0   -{\frac {8}{9}}\,{\mu}^{2} K_1^{4}+\frac{10}{3}\,{\mu}^{4}   K_1 ^{2}+10\,  K'_0 {\mu}^{4}-{\frac{28}{27}}\,{\mu}^{3}K_1^{3}-
8\,{\mu}^{2}  {K'_0}^{2}+{\frac {16}{9}}\,   K_0    ^{2}{K'_1}^{2}+16\,  K_0^{4}  K_3   ^{2}+{\frac {16}{9}}\,  K_1 ^{2} {(K'_0)}^{2}-14\,{\mu}^{3}  K_0  K_3   -6\,{\mu}^{2}
K_0^{2} K'_2   +\frac{14}{3}\,{\mu}^{3} K_0   K'_1.  
$

Though the formula for the determinant  looks   ugly, for explicit $K_i$ given at the beginning of this paragraph,  one   can   calculate it (Maple$^\registered$ does it within few  seconds, a human needs around one hour for it). 
Calculating this  formula for the $K_i$ corresponding to the projective 
connections corresponding to the metrics \ref{1a}, \ref{1b}, \ref{1c} 
(explicit formulas for the projective connections are at the beginning of the present paragraph), we obtain that the result   is not  zero implying the system \eqref{lin1} corresponding to the projective connections corresponding to the metrics $g$  from cases 1a-- 2c of theorem 1 does not admit a three dimensional $\mathcal{A}$ under the additional assumption that  the matrix of  $L_v$ is
the first  matrix 
of \eqref{mat3bis}.  

\subsection{ Calculations related to the proof of Theorem 2
  for the  metric 1a    from Theorem 1  assuming that the matrix of $L_v$ is as the second matrix of   
 (51) } \label{case3} 

   For the metric \ref{1a}, 
 we consider the coordinate system  
 $x_{new}= \frac{x_{old} + y_{old}}{2}$, $y_{new}=\frac{ x_{old}- y_{old}}{2}$. 

 In this coordinate system, 
 the projective  vector field is $\frac{\partial }{\partial x}$, and the projective connection of $g$ is \eqref{pcb}. 
  Direct calculations shows that the matrix of $\bar a$ is given by  
      \begin{equation} \label{bara} \left( \begin {array}{cc} 
      {\frac {{e^{4\,x}} \left( {c}\,
{e^{-3\,x-3\,y}}(y-x)-{e^{-3\,x+3\,y}}(x+y)\right) }{4\sqrt [3]{y}{{c}}^{\frac{2}{3}}}}&
{\frac {{
e^{4\,x}} \left( {c}\,{e^{-3\,x-3\,y}}(y-x)+
{e^{-3\,x+3\,y}}(x+y) \right) }{4\sqrt [3]{y}{{c}}^{\frac{2}{3}}}}\\\noalign{\medskip}{\frac {{
e^{4\,x}} \left( {c}\,{e^{-3\,x-3\,y}}(y-x)+
{e^{-3\,x+3\,y}}(x+y) \right) }{4\sqrt [3]{y}{{c}}^{\frac{2}{3}}}}&{\frac {{e^{4\,x}}
 \left( {c}\,{e^{-3\,x-3
\,y}}(y-x)-
{e^{-3\,x+3\,y}}(x+ y)
 \right) }{4\sqrt [3]{y}{{c}}^{\frac{2}{3}}}}\end {array} \right).\end{equation} 
 Direct computations shows that the following matrix  $P$ is a partial solution of the equation \eqref{cb}. 
\begin{equation} \label{P} 
  \left( \begin {array}{cc} {\frac {x \left( -2\,y{e^{3\,y}}-c\,x{e^{-3\,y}}+2\,yc\,{e^{-3\,y}}-x{e^{3\,y}} \right) {e^{x}}}{8\sqrt [3]{y}{c}^{\frac{2}{3}}}}&{\frac {x \left( 2\,y{e^{3\,y}}-c\,x{e^{-3\,y}}+2\,yc\,{e^{-3\,y}}+x{e^{3\,y}} \right) {e^{x}}}{8\sqrt [3]{y}{c}^{\frac{2}{3}}}}\\\noalign{\medskip}{\frac {x \left( 2\,y{e^{3\,y}}-c\,x{e^{-3\,y}}+2\,yc\,{e^{-3\,y}}+x{e^{3\,y}} \right) {e^{x}}}{8\sqrt [3]{y}{c}^{\frac{2}{3}}}}&{\frac {x \left( -2\,y{e^{3\,y}}-c\,x{e^{-3\,y}}+2\,yc\,{e^{-3\,y}}-x{e^{3\,y}} \right) {e^{x}}}{8\sqrt [3]{y}{c}^{\frac{2}{3}}}}\end {array} \right) .
\end{equation}
Thus,  the solution of the equation \eqref{cb} has the form \eqref{ansatz}. 
 Substituting this ansatz in the  equations \eqref{lin1} and solving  the last three equations with respect to the first derivatives, we obtain the system

\begin{equation}  \label{s2} \left. 
\begin{array}{ccccc} 0&=&
{ a_{11}}  -\frac{2}{3}\,{ K_1}  { a_{11}}  +2\,{ K_0}  { a_{12}} &+& {\frac {{y{{}}}^{\frac{2}{3}} \left( -{e^{3y }}+{ c}\,{e^{-3
\,y}} \right) }{4{{ c}}^{\frac{2}{3}}}} \\
 { a'_{11}}  & = &  - 2\,{ a_{12}}
  +\frac{4}{3}\,{ K_2} 
   { a_{11}}  -\frac{2}{3}\,{ K_1} { a_{12}}  -2\,{ K_0} { a_{22}} &-&  {\frac {{y{{}}}^{\frac{2}{3}} \left( {e^{3y }}+{ c}\,{e^{-3
\,y}} \right) }{4{{ c}}^{\frac{2}{3}}}} \\
{ a'_{12}} & =&  -\tfrac{1}{2}{ a_{22}}+{ K_3}  { a_{11}}  +\frac{1}{3}\,{ K_2}  { a_{12}} -\frac{2}{3}\,{ K_1} 
{ a_{22}} &+&  {\frac {{y{{}}}^{\frac{2}{3}} \left( {e^{3y }}-{ c}\,{e^{-3
\,y}} \right) }{4{{ c}}^{\frac{2}{3}}}}\\  
 { a'_{22}}& =&  2\,{ K_3} { a_{12}} -\frac{2}{3}\,{ K_2}  
 { a_{22}}, 
\end{array}\right\}\end{equation}
 where $K_i$ are the coefficients of the projective connection 
  \eqref{pcb}. 

  The first equation of \eqref{s2} plays the role of the first equation of \eqref{tmp1}. 
  Differentiating the first equation of \eqref{s2} by  $y$ and substituting the values of  derivatives from the last three equations of  \eqref{s2}  inside, we obtain the following   nonhomogeneous linear equation in $a_{ij}$, which plays  the role of the second  equation of \eqref{tmp1}. \\
  {\tiny 
$
-\left( 38\,\sqrt [3]{y}{{c}}^{\frac{5}{3}}{e^{18\,y}}-108\,
{y}^{4/3}{{c}}^{\frac{11}{3}}{e^{6\,y}}+72\,{y}^{4/3}{{c}}^{\frac{8}{3}}{e^{12\,y}}-38\,\sqrt [3]{y}{{c}}^{\frac{11}{3}
}{e^{6\,y}}-108\,{y}^{4/3}{{c}}^{\frac{5}{3}}{e^{18\,y}}-9\,\sqrt [3]{y}{{c}}^{\frac{2}{3}}{e^{24\,y}}+9\,{{c}}^{14/3}\sqrt [3]{y} \right) {a_{11}} $ \\$ - \left( 72\,{y}^{4/3}{{c}}^{\frac{11}{3}}{e^{6\,y}}+20\,\sqrt [3]{y}{{c}}^{\frac{5}{3}}{e^{18\,
y}}-72\,{y}^{4/3}{{c}}^{\frac{5}{3}}{e^{18\,y}}-9\,
\sqrt [3]{y}{{c}}^{\frac{2}{3}}{e^{24\,y}}+20\,\sqrt [3]{y}{{c}}^{\frac{11}{3}}{e^{6\,y}}-9\,{{c}}^{14/3}\sqrt [3]{y}+58\,\sqrt [3]{y}{{c}}^{\frac{8}{3}}{e^{12\,y}} \right){a_{12}} \\ - \left( 72\,{y}^{4/3}{{c}}^{\frac{8}{3}}{e^{12\,y}}
+36\,{y}^{4/3}{{c}}^{\frac{5}{3}}{e^{18\,y}}+36\,{y}^{4/3}{{c}}^{\frac{11}{3}}{e^{6\,y}} \right){a_{22}} \\ =    -25\,{{c}}^{3}{e^{9\,y}}{y}^{2}+72\,{{c}}^{3}{y}^{3}{e^{9\,y}}+72\,{y}^{3}{e^{15\,y}}{{c}}^{2}+9\,{e^{21\,y}}{y}^{2}{c}+25\,{e^{15\,y}}{y}^{2}{{c}}^{2}-9\,{{c}}^{4}{e^{3\,y}}{y}^{2}  .
$} 

Differentiating this equation by $y$ and  substituting the values of  derivatives from the last three equations  of \eqref{s2}  inside, we obtain the following   nonhomogeneous linear equation on $a_{ij}$, which plays  role of the third  equation of \eqref{tmp1}. \\ {\tiny 
$
 \left( -132\,{e^{24\,y}}{{c}}^{\frac{5}{3}}+103\,{e^{18\,y}}
{{c}}^{\frac{8}{3}}+132\,{{c}}^{{\frac {17}{3}}}+27\,{e^{30\,y}}{{c}}^{\frac{2}{3}}+2640\,{{c}}^{\frac{14}{3}}y{e^{6\,y}}+252\,{{c}}^{{\frac {17}{3}}}y+1128\,{e^{12\,y}}y{{c}}^{\frac{11}{3}}-103\,{{c}}^{\frac{14}{3}}{e^{6\,y}}\right. $ \\ $\left. -27\,{{c}}^{{\frac {20}{3}}}{e^{-6\,y}} +6048\,{y}^{2}{{c}}^{\frac{14}{3}}{e^{6\,y}}-6912\,{y}^{2}{{c}}^{\frac{11}{3}}{e^{12\,y}}+1980\,{e^{24\,y}}y{{c}}^{\frac{5}{3}}+14688\,{e^{18\,y}}{y}^{2}{{c}}^{\frac{8}{3}}-4656\,{{c}}^{\frac{8}{3}}{e^{18\,y}}y \right) {a_{11}} \\  + \left( -144\,{{c}}^{{\frac{17}{3}}}y+8640\,{e^{18\,y}}{y}^{2}{{c}}^{\frac{8}{3}}-78\,{{c}}^{{\frac {17}{3}}}-2592\,{{c}}^{\frac{14}{3}}y{e^{6\,y}}-107\,{{c}}^{\frac{14}{3}}{e^{6\,y}}-107\,{e^{18\,y}}{{c}}^{\frac{8}{3}}+3456\,{y}^{2}{{c}}^{\frac{11}{3}}{e^{12\,y}} +27\,{e^{30\,y}}{{c}}^{\frac{2}{3}}\right. $ \\ $\left. -5568\,{e^{12\,y}}y{{c}}^{\frac{11}{3}}+27\,{{c}}^{{\frac {20}{3}}}{e^{-6\,y}}-5184\,{y}^{2}{{c}}^{\frac{14}{3}}{e^{6\,y}}+1872\,{e^{24\,y}}y{{c}}^{\frac{5}{3}}-4\,{e^{12\,y}}{{c}}^{\frac{11}{3}}-1248\,{{c}}^{\frac{8}{3}}{e^{18\,y}}y-78\,{e^{24\,y}}{{c}}^{\frac{5}{3}} \right) {a_{12}} \\ - \left( 6048\,{e^{18\,y}}{y}^{2}{{c}}^{\frac{8}{3}}+336\,{{c}}^{\frac{8}{3}}{e^{18\,y}}y+6912\,{y}^{2}{{c}}^{\frac{11}{3}}{e^{12\,y}}+108\,{e^{24\,y}}y{{c}}^{\frac{5}{3}}+864\,{y}^{2}{{c}}^{\frac{14}{3}}{e^{6\,y}}+456\,{e^{12\,y}}y{{c}}^{\frac{11}{3}}+108\,{{c}}^{{\frac {17}{3}}}y+336\,{{c}}^{\frac{14}{3}}y{e^{6\,y}} \right) {a_{22}} \\ =  1872\,{e^{21\,y}}{y}^{\frac{8}{3}}{{c}}^{2}+57\,{e^{21\,y}}{y}^{\frac{5}{3}}{{c}}^{2}- 27\,{{c}}^{6}{y}^{\frac{5}{3}}{e^{-3\,y}}+27\,{e^{27\,y}}{y}^{\frac{5}{3}}{c}+286\,{{c}}^{3}{e^{15\,y}}{y}^{\frac{5}{3}}+286\,{{c}}^{4}{e^{9\,y}}{y}^{\frac{5}{3}}+5184\,{{c}}^{4}{e^{9\,y}}{y}^{\frac{11}{3}}-57\,{{c}}^{5}{y}^{\frac{5}{3}}{e^{3\,y}}\\ +8640\,{e^{15\,y}}{y}^{\frac{11}{3}}{{c}}^{3}+4944\,{e^{15\,y}}{y}^{\frac{8}{3}}{{c}}^{3}+144\,{{c}}^{4}{e^{9\,y}}{y}^{\frac{8}{3}}144\,{{c}}^{5}{y}^{\frac{8}{3}}{e^{3\,y}}
$}\\
Direct calculations show that   det$\left(
\begin{array}{ccc}  b_1  &  m_{12} &    m_{13}  \\ 
b_{2} & m_{22}   &    m_{23}   \\
 b_{3} & m_{32} &   m_{33}   \end{array} \right)   $   is equal to  \\
\begin{equation}\label{ror}
\textrm{\hspace{-.2cm}\tiny${ {9}}\,{\frac {18\,{e^{-9\,y}}{{c}}^{4}y-18
\,{{c}}^{2}{e^{3\,y}}y+32\,{{c}}^{3}{e^{-3\,y{{
}}}}+32\,{{c}}^{2}{e^{3\,y}}-18\,{e^{9\,y}}y{c}+9\,{e^{-9\,y}}{{c}}^{4}+18\,{{c}}^{3}{e^{-3\,y}}y-{{c}}^{5}{e^{-15\,y}}+
9\,{e^{9\,y}}{c}-{e^{15\,y}}}{{64}{{c}}^{\frac{8}{3}}{y}^{7/3}}}.$} 
\end{equation}

  We see that it  is not zero, which gives us a contradiction which proves Theorem \ref{main2} for the metric from the case  \ref{1a} of Theorem \ref{main}.

\subsection{ Proof of Theorem \ref{main2}  for the metrics  \ref{3a}, \ref{3b}, \ref{3c},   and \ref{3d} } \label{case1}

In all these cases the metric has the form $2(Y(y)+x)dxdy $. We first explain that  the dimension of   the space $\mathcal{A}$ coincides with the dimension of  the space of integrals quadratic in momenta for the Hamiltonian 
$$H:T^*D\to \mathbb{R} \ , \ \  \  H(x,y, p_x, p_y):= \frac{p_xp_y}{Y(y)+x}.$$

Indeed, as we explained in \S\ref{2.1}, for every solution $a$  
the function $ (\det(g))^{2/3}  a_{ij}\xi^i \xi^j $ is an integral of the geodesic flow of $g$, and vice versa.  
Since  the mapping $a\mapsto  (\det(g))^{2/3}  a$ is  linear and bijective, the dimensions of the space $\mathcal{A}$ and of    the space of integrals quadratic in momenta coincide. 

Note that the space of the integrals quadratic in momenta  is at least two-dimensional. Indeed, 
every linear combination of the Hamiltonian $H$ and of the integral coming from the projectively equivalent metric \eqref{E3} by formula \eqref{integral}  is  an integral. In the notations below, these integrals will correspond to $a= \const$, $c=0$. Our goal is to show that in the cases \ref{3a}, \ref{3b}, \ref{3c} all integrals have $a= \const$, $c=0$, and that in the case \ref{3d} there exists an additional linearly independent  integral.

Suppose a function $f:T^*D\to \mathbb{R}$ of the form $a(x,y)p_x^2+ b(x,y) p_xp_y +  c(x,y) p_y^2$ 
 is an integral for the geodesic flow of $g$. Then, the condition 
 $0=\{H,f\}$, after multiplication by $-(Y+x)^2$,   reads \begin{eqnarray*}
0&=& -(Y+x)^2\cdot \left\{\frac{p_xp_y}{Y+x}, ap_x^2+ bp_xp_y+ cp_y^2\right\}  \\ &=& p_x^3(Y+x)a_y + p_x^2 p_y ((Y+x)a_x + (Y+x)b_y + 2  a + Y' b)\\ &+& p_xp_y^2 ((Y+x)b_x + (Y+x)c_y+  b + 2 Y' c )+ p_y^3 (Y+x)c_x,\end{eqnarray*} i.e., is equivalent to the following system of PDE: 
\begin{equation}\label{syst} 
 \left.\begin{array}{rcc} a_y&=&0 \\ (Y+x)a_x + (Y+x)b_y + 2 a + Y' b&=&0\\  
 (Y+x)b_x + (Y+x)c_y+  b + 2 Y'c &=&0\\ c_x&=&0 \end{array}
 \right\}\end{equation} 
We see that the first (the last, respectively)  equation of \eqref{syst}  implies  that the 
function  $a$ ($c$, respectively) is a function of the variable $x$ ($y$, respectively) only.  

Solving the second and the third equations with respect to the derivatives of $b$, we obtain 
$$
  b_y=-\frac{2a + Y'b}{Y+x}- a' \  , \  \  \  \  \  
  b_x=-\frac{b + 2 Y'c}{Y+x}- c'.
$$ 
 Substituting these  
 expressions for the derivatives of $b$ in the identity $\frac{\partial b_x}{\partial y}-  
  \frac{\partial b_y}{\partial x}=0$, we obtain 
   \begin{equation}\label{sys2}  - a'' 
  x+ c''   x+ c''   Y -3a' +3c'  Y' - a''   Y +2 Y''  c =0. \end{equation} Taking the $\frac{\partial^3}{\partial^2 x\partial y}-$derivative 
   of this equation, we obtain ${a''''}\,  Y'=0$. Since the function $Y$  from cases \ref{3a}, \ref{3b},  \ref{3c}, \ref{3d}  is not constant, we can assume 
    $Y'\ne 0$.  Then,  $a= \alpha_3 x^3+ \alpha_2 x^2 + \alpha_1 x  + \alpha_0$, where $\alpha_i$ are constants.
    Substituting this in \eqref{sys2}, we obtain 
 \begin{equation} \label{sys4} -15\alpha_{3}x^2+ \left(
c'' -8\alpha_{2}-6Y 
\alpha_{3} \right) x-3\alpha_{1}+3c' Y' + c''
   Y  +2 Y''   c -2Y \alpha_{{2}}=0.\end{equation}
The left-hand side of this equation   is a polynomial in $x$ whose coefficients depend on $y$ only. They must be zero implying $\alpha_3=0$, $c=4\alpha_2y^2 + \beta_1y + \beta_0$. Then, the equation \eqref{sys4} reads 
 \begin{equation} \label{sys5} 
 6Y  \alpha_{2}-3\alpha_{1}+\left( 3\beta_{1}+24\alpha_{2}y \right)Y' +(2 \beta_{0}+2
 \beta_{1}y+8     \alpha_{{2}}{y}^{2})Y'' =0. 
\end{equation}
If $\alpha_{1}= \alpha_{2}=\beta_{1}=\beta_0= 0$, then    $c=0$, $a= \const$ implying that the integral  is a linear combination of the Hamiltonian and   the integral coming from the projectively  equivalent metric \eqref{E3} by formula \eqref{integral}.   Otherwise 
\eqref{sys5}  is an ODE for the function $Y$. 
Substituting the functions $Y$ from the cases 
  \ref{3a}, \ref{3b}, \ref{3c}  of   Theorem \ref{main}  we see that they are not solutions of this ODE. Thus, the metrics from the cases 
  \ref{3a}, \ref{3b}, \ref{3c}  from  Theorem \ref{main} have  2-dimensional $\mathcal{A}$. Substituting the functions $Y$ from the case 
   \ref{3d}  from  Theorem \ref{main}, 
     we see that it is a solution of this ODE if and only if $\beta_1= \alpha_2=0$, $4 \beta_0 = 3 \alpha_1$.  We see that there is precisely one additional parameter ($\beta_0$) we can freely choose to construct  the integral,
      i.e., the space of the integrals  is at most three-dimensional.  Direct calculations show that as the additional    integral  we  can take the integral corresponding to the metric $\tilde g$ by formula \eqref{integral}. 

\section{ Proof of Theorem \ref{main3}} \label{pr4} 

The goal is to show that no metric from Theorem \ref{main}
has a Killing vector field. 
 It is sufficient to do it for the cases \ref{1a} -- \ref{2c} only, since  in view of 
  \S\ref{case1}, in  the cases
  \ref{3a} -- \ref{3d}  we know the space of quadratic integrals  of the metrics \ref{3a} -- \ref{3d}, so it is  sufficient to check  that  no quadratic integral    is degenerate at every point, which is an easy exercise.  Moreover, in the case \ref{3d} the space of  quadratic integrals  is precisely 
  3-dimensional implying the metric  admits no Killing vector field, see  \cite[Section 5]{kruglikov}. 

   We will use the following approach which was known to Darboux  \cite[\S\S688,689]{Darboux} and Eisenhart \cite[pp. 323--325]{eisenhartsurfaces}, 
   see also \cite{kruglikov} for an equivalent approach leading to  similar calculations.  
 For every  $g$ from the cases \ref{1a} -- \ref{2c} of Theorem \ref{main},  let us  
consider  the  following functions on $D^2$:  \begin{itemize}  
\item The scalar curvature $R:=\sum_{i,j,k} R^i_{ijk} g^{jk},$ where $R^i_{h jk }$ is the curvature tensor of $g$,  
\item The  square of the length of the derivative of the scalar curvature
    $L:=\sum_{i,j} g^{ij}
    \frac{\p R}{\p x_i}\frac{\p R}{\p x_j}$,
\item The laplacian of the scalar curvature $\Delta
     :=\frac{1}{\sqrt{\det(g)}} \sum_{i,j} \frac{\p }{\p x_i}
     \left( g^{ij} \sqrt{\det(g)}\frac{\p R}{\p x_j}\right).$
\end{itemize}  
If the metric admits a Killing  vector field $K$, 
then in a coordinate system  $(x_1,x_2)$  such that $K= \frac{\partial }{\partial x_1}$, 
all these functions depend of $x_2$ only. Then, the differentials $dR$, $dL$  are proportional, and the differentials $dR$, $d\Delta$ are proportional. Then, in every coordinate system $(x,y)$  the following 
     determinants  are zero: 

     \begin{equation} 
    \det\begin{pmatrix} \frac{ \partial  R}{\partial x}  & \frac{ \partial  R}{\partial y} \\ \frac{ \partial  L}{\partial x}  & \frac{ \partial  L}{\partial y} \end{pmatrix}\  , \  \   \det\begin{pmatrix} \frac{ \partial  R}{\partial x}  & \frac{ \partial  R}{\partial y} \\ \frac{ \partial  \Delta}{\partial x}  & \frac{ \partial  \Delta}{\partial y} \end{pmatrix}  .\end{equation} 

     Calculating these determinants for all metrics from the cases \ref{1a} -- \ref{2c}, we see that  in every case they are   not
      zero  implying that the metrics admit no Killing vector field. Note that it is sufficient to calculate the determinants for the metrics from the  cases \ref{1a} -- \ref{1c} only, since the cases  
      \ref{2a} -- \ref{2c}, up to multiplication by a constant, can be obtained from  the cases \ref{1a} -- \ref{1c} by replacing $x$ by $z$ and $ y$ by $\bar z$.

\begin{appendix}

\section{  Appendix: Dini's theorem for pseudo-Riemannian metrics}

\centerline{\bf  by 
Alexei V. Bolsinov\footnote{
Department of Mathematical Sciences,  
Loughborough University,  LE11 3TU 
UK, A.Bolsinov@lboro.ac.uk}, 
Vladimir S. Matveev\footnote{ Institute of Mathematics, FSU Jena, 07737 Jena Germany,  vladimir.matveev@uni-jena.de}, and  
Giuseppe Pucacco\footnote{Dipartimento di Fisica Universit\`a di Roma ``Tor Vergata", 00133 Rome Italy, pucacco@roma2.infn.it}}

\vspace{2ex} 

\subsection{Introduction} 

Consider a  Riemannian or a  pseudo-Riemannian metric $g=(g_{ij})$ on  a surface $M^2$.  We say that a  metric $\bar g$ on the same surface is \emph{projectively equivalent} to $g$,
 if every geodesic of $\bar g$ is a reparameterized geodesic of $g$.   
In 1865 Beltrami  \cite{Beltrami} asked\footnote{ Italian original from \cite{Beltrami}: 
La seconda $\dots$  generalizzazione $\dots$ del nostro problema,   vale a dire:   riportare i punti di una superficie sopra un'altra superficie in modo  che alle linee geodetiche della prima corrispondano linee geodetiche della seconda.}  to describe all pairs of  projectively equivalent Riemannian  metrics on  surfaces. 
From the context it is clear that he considered this problem locally, in a neighborhood of almost every  point.  

Theorem A below, which is the main result of this note (which is a short version of \cite{pucacco}),  gives  an answer to the following generalization of  the question of Beltrami: we allow the metrics $g$ and $\bar g$ to be pseudo-Riemannian.

\begin{Tha} \label{thm3}
Let $g$, $\bar g$ be projectively equivalent  metrics on $M^2$, and $\bar g\ne \textrm{const} \cdot g$ for every $\textrm{const}\in \mathbb{R}$. Then, in the neighborhood of almost every point there exist coordinates $(x,y)$  such that the metrics are  as in the following table. \\
 \begin{tabular}{|c||c|c|c|}\hline &  \textrm{Liouville Case} & \textrm{Complex-Liouville case} & \textrm{Jordan-block case}\\ \hline \hline
$g$ & $(X(x)-Y(y))(dx^2 \pm  dy^2)$ &  $2\Im(h)dxdy$ & $\left( 1+{x} Y'(y)\right)dxdy $
\\  \hline  $ \bar g$  &$ \left( \frac{1}{Y(y)}-\frac{1}{X(x)}\right) \left( \frac{dx^2}{X(x)} \pm   \frac{dy^2}{Y(y)} \right)$& 
 \begin{minipage}{.3\textwidth}$-\left(\frac{\Im(h)}{\Im(h)^2 +\Re(h)^2}\right)^2dx^2 \\   +2\frac{\Re(h) \Im(h)}{   (\Im(h)^2 +\Re(h)^2)^2} dx dy  \\ +  \left(\frac{\Im(h)}{\Im(h)^2 +\Re(h)^2}\right)^2dy^2 $
 \end{minipage} &  \begin{minipage}{.3\textwidth}$  \frac{1+{x} Y'(y)}{Y(y)^4} \bigl(- 2Y(y) dxdy\\
    + (1+{x} Y'(y))dy^2\bigr)$\end{minipage}\\ \hline 
\end{tabular} 
where $h:=\Re(h) + i\cdot \Im(h)$ is  a holomorphic function of the variable $z:= x +i\cdot y $.   \end{Tha} 

\begin{rema} 
It it natural to consider the metrics from the Complex-Liouville case as the complexification of the metrics from the Liouville case: indeed, in the complex coordinates $z= x + i\cdot y$, $\bar z= x - i\cdot y$, the metrics  have the form

 \begin{equation*}  \left. \begin{array}{ccc}   ds^2_g & =  &  -\tfrac{1}{4}( \overline{h( z}) - h(z) )\left(d\bar z^2 -   dz^2\right),\\   
     ds^2_{\bar g} & =  &  -\tfrac{1}{4}\left(\frac{1}{\overline{h( z})} - \frac{1}{h(z)} \right)\left(\frac{d\bar z^2}{ \overline{h(z})} - \frac{ dz^2}{h(z)}\right)\end{array}\right. \end{equation*}(this form is used in the proof of Theorem \ref{main}). 
\end{rema} 

\begin{remb}  In the Jordan-block case, 
if $dY\ne 0$ (which is always the case  at almost every point,  if the restriction of  $g$ to any neighborhood   does   not admit a Killing vector field), after a  local coordinate change, the metrics  $g$ and   $\bar g$  have the form   \begin{eqnarray*}ds_g^2 & =&  \left( \tilde Y(y) +{x}{} \right)dxdy \\ 
   ds_{\bar g}^2   & =&  -\frac{2(\tilde Y(y)+x)}{y^3}dxdy  + \frac{(\tilde Y(y)+x)^2}{y^4}dy^2
 \end{eqnarray*}
 (this form is used in the proof of Theorem \ref{main}). 
  \end{remb} 

We see that  the metric $g$ 
 from Complex-Liouville and  Jordan-block cases  always has signature (+,--), and the metric $g$ 
 from  the Liouville case has  signature $(+,+)$ or $(-,-)$, if the sign ``$\pm$" is ``$+$".  
In this case, 
the formulas from Theorem A  are precisely the formulas obtained by Dini in \cite{Dini}.

We do not insist that we are the first to find these normal forms of projectively equivalent pseudo-Riemannian  metrics. According to \cite{Aminova},   a description of projectively equivalent metrics  was obtained by 
  P. Shirokov in \cite{shirokov}. Unfortunately, we were not able to find the reference \cite{shirokov} to check it.  The result of Theorem A could be     even more classical, see  Remark D.

Given  two projectively  equivalent metrics,  it is easy to understand  what case they belong to. Indeed, the $(1, 1)$-tensor 
  $
  G^i_{j}:= \sum_{\alpha=1}^{2} \bar g_{j\alpha} g^{i\alpha}, 
  $   where $ g^{i\alpha}$ is inverse to $ g_{i \alpha}$,  
   has two different real eigenvalues in the Liouville  case, two complex-conjugate eigenvalues in the Complex-Liouville  case, and is (conjugate to) a Jordan-block in the Jordan-block case.

 There exists  an interesting and useful  connection of projectively equivalent metrics with integrable systems. 

Recall that 
 a function $F:{T^*}M^2\to \mathbb{R}$ is called \emph{an integral} of the geodesic flow of $g$, if $\{ H, F\}=0$, where $H:= \tfrac{1}{2} g^{ij} p_ip_j:T^*M^2\to \mathbb{R}$ is the kinetic energy corresponding to the metric, and $\{ \ , \ \}$ is the standard Poisson bracket on $T^*M^2$. Geometrically, this condition means 
that the function is constant on the orbits of the Hamiltonian system with the Hamiltonian $H$.  We say the  integral $F$ is {\it quadratic in  momenta,} if for  every local coordinate system $(x,y)$ on $M^2$  it has the form 
\begin{equation} \label{integral1} F(x,y, p_x, p_y)= 
a(x,y)p_x^2+ b(x,y)p_xp_y+ c(x,y)p_y^2  
\end{equation} in the canonical coordinates $(x,y,p_x, p_y)$ on $T^*M^2$. 
Geometrically, the  formula \eqref{integral1} 
means that the restriction of the integral to every cotangent space $T^*_pM^2\equiv \mathbb{R}^2$  is a homogeneous quadratic function. Of course, $H$ itself  is  an integral quadratic  in the momenta for $g$. We will say that the integral $F$ is {\it nontrivial}, if $F\ne  \textrm{const} \cdot H$ for all $\textrm{const} \in \mathbb{R}$. 

\begin{Thb} \label{mainb}  
 Suppose the metric $g$ on  ${M}^2$ admits a nontrivial integral quadratic in momenta. Then, in a neighbourhood of almost every point there exist coordinates  $(x,y)$ such that  the metric and the integral are as in the following table
\begin{center}\begin{tabular}{|c||c|c|c|}\hline & \textrm{Liouville Case} & \textrm{Complex-Liouville Case} & \textrm{Jordan-block Case}\\ \hline \hline 
$g$ & $(X(x)-Y(y))(dx^2 \pm dy^2)$ &  $\Im(h)dxdy$ & $\left(1+{x} Y'(y)\right)dxdy $
\\ \hline  $F$&$\tfrac{X(x)p_y^2 \pm Y(y)p_{x}^2}{X(x)-Y(y)} $&  
 $p_x^2 - p_y^2   + 2\tfrac{\Re(h)}{\Im(h)}p_xp_y, $ & $p_x^2 -2  \frac{Y(y)}{1+ {x}{} Y'(y) }p_xp_y$  \\ \hline     
\end{tabular} \end{center}
$\textrm{where $h:=\Re(h) + i\cdot \Im(h)$ is a}$  \textrm{holomorphic function} $\textrm{of  the variable $z:= x+ iy$.}$
\end{Thb}

  Indeed, as was shown in  \cite{MT,dim2},  and as it was essentially known to  Darboux~\cite[\S\S600--608]{Darboux}, if two metrics $g$ and $\bar g$ are projectively equivalent, then \begin{equation}\label{Integral}  I:TM^2\to \mathbb{R}, \
\
 I(\xi):=\bar g(\xi,\xi)\left(\frac{\det(g)}{\det(\bar g)}\right)^{2/3}
\end{equation}
 is an integral of the geodesic flow of  $g$. Moreover, it was shown in \cite[\S 2.4]{bryant},  see also \cite{quantum},    the above statement  is proved to be  true\footnote{with a good will, one also can attribute this result to Darboux} in the  other direction: if the function \eqref{integral1} is an integral for the geodesic flow of $g$,  then the metrics $g$ and $\bar g$ are projectively equivalent. Thus, Theorem A and Theorem B are equivalent. 
In this paper,  we will actually prove Theorem B obtaining Theorem A  as its consequence.

 \begin{remc} The corresponding natural Hamiltonian problem on the hyperbolic plane has  been recently treated in \cite{PR} following the approach used by Rosquist and Uggla \cite{ru:kt}. \end{remc} 

 \begin{remd}  The formulas that will appear in the proof  are very close to those in \S 593 of \cite{Darboux}. Darboux  worked over the complex numbers and therefore did  not care about whether the metrics are Riemannian or  pseudo-Riemannian. For example, the Liouville and Complex-Liouville case are the same for him.  Moreover,   in \S 594, Darboux gets  
the  formulas that are very close to that of the Jordan-block case, though he was interested  in the Riemannian case only, and, hence, treated this ``imaginary" case as not interesting.  \end{remd}

\subsection{Proof of Theorem B  (and, hence, of Theorem A)} 
If the metric $g$ has signature (+, +) or (--,--), Theorem A and, hence, Theorem B,  were obtained by Dini in \cite{Dini}. Below we assume that the metric $g$ has signature (+,--).

\subsubsection{Admissible coordinate systems and 
Birkhoff-Kolokoltsov forms} \label{admissible}

Let $g$ be a pseudo-Riemannian metric on $M^2$  of signature $(+, -)$. 
Consider (and fix)  two linearly independent  vector fields $V_1, V_2$ on $M^2$  such that 
\begin{itemize}  
 \item  $g(V_1, V_1) =g(V_2, V_2)=0$ and 
 \item $g(V_1, V_2)>0$.
 \end{itemize} 
Such vector fields always exist locally (and, since our result is local, this is sufficient for our proof).

We will say that a local  coordinate system $(x,y)$
 is {\it admissible}, if  the vector fields  $\frac{\partial }{\partial x}$ and  $\frac{\partial }{\partial y} $ are proportional to $V_1, V_2$ with positive coefficient of proportionality:  $$V_1(x,y)=  \lambda_1(x,y)\frac{\partial }{\partial x} , \  \  \ V_2(x,y) =  \lambda_2(x,y) \frac{\partial }{\partial y}, \  \  \ \textrm{where $\lambda_i>0, \;\; i=1,2$}.$$ 
 Obviously, 
\begin{itemize} 
\item  admissible coordinates exist in a sufficiently small neighborhood  of every point, \item the metric $g$  in  admissible coordinates has the form 
 \begin{equation}\label{metric} 
 ds^2 =f(x,y)dxdy , \  \  \ \textrm{where $f>0$},  \end{equation}       
     \item two admissible  coordinate systems  in 
      one neighborhood  are connected by \begin{equation} \label{coordinatechange} \begin{pmatrix} x_{new}\\
 y_{new}\end{pmatrix} 
 = \begin{pmatrix} x_{new}(x_{old}) \\
 y_{new}(y_{old})\end{pmatrix} , \  \ \textrm{where  $\frac{dx_{ new}}{dx_{old}}>0$, $\frac{dy_{ new}}{dy_{old}}>0$}. \end{equation}   
\end{itemize}

\begin{LemmaA}  \label{BK} 
Let $(x,y)$ be an admissible coordinate system for $g$. 
Let  $F$ given by  \eqref{integral1} be  an  integral for  $g$. 
Then,    $B_1:= \frac{1}{\sqrt{|a(x,y)|}}dx$  ($B_2:= \frac{1}{\sqrt{|c(x,y)|}}dy$, respectively) is a    1-form, which is defined at points such that $a\ne 0$ ($c\ne 0$, respectively).  Moreover, the coefficient 
 $a$ ($c$, respectively) depends only on $x$ ($y$, respectively), which in particular imply that the forms $B_1$, $B_2$ are  closed.   \end{LemmaA} 

\begin{reme}  The forms $B_1, B_2$ are not the direct analog of the ``Birkhoff" 2-form introduced by Kolokoltsov in \cite{Kol}. In a certain sense, they are  the real analogue of the different branches of the  square root of the Birkhoff  form.   
\end{reme} 
{\bf Proof of  Lemma A.} The first part of the statement, namely that the $\frac{1}{\sqrt{|a|}} dx$ ($\frac{1}{\sqrt{|c|}}dy$, respectively) transforms as a $1$-form under admissible coordinate changes is evident:  indeed, after the coordinate change  
\eqref{coordinatechange}, the  momenta transform as follows: 
 $p_{x_{old}}= p_{x_{new}}\frac{d{x_{new}}}{d{x_{old}}}$, $p_{x_{old}}= p_{x_{new}}\frac{d{x_{new}}}{d{x_{old}}}$. Then, the integral $F$ in the new coordinates has 
  the form $$ \underbrace{{a}\left(\frac{d{x_{new}}}{d{x_{old}}}\right)^2}_{a_{new} } {p_{x_{new}}^2}  + \underbrace{{b}\frac{d{x_{new}}}{d{x_{old}}}\frac{d{y_{new}}} {d{y_{old}}}}_{b_{new}}  {p_{x_{new}}} {p_{y_{new}}} + \underbrace{{c}\left(\frac{d{y_{new}}}{d_{y_{old}}} \right)^2}_{c_{new}} {p_{y_{new}}^2}.$$  Then, the formal expression $\frac{1}{\sqrt{|a|}}dx_{old}$
  $ \left(  \frac{1}{\sqrt{|c|}}dy_{old}, \textrm{ respectively}\right)$   transforms in  $$\frac{1}{\sqrt{|a|}} \frac{d{x_{old}}}{d{x_{new}}} dx_{new} \ \ \ \ \  \left(\textrm{  $ \frac{1}{\sqrt{|c|}}\frac{d{y_{old}}}{d{y_{new}}}dy_{new}$, respectively}\right),  $$ 
  which is precisely the transformation law of  1-forms.

Let us prove that the forms are closed. 
 If $g$ is given by \eqref{metric}, its Hamiltonian $H$ is given by $\frac{p_xp_y}{2f}$, and the condition $0=\{H, F\}$ reads \\
 $\begin{array}{ccl}
0&= &\left\{\frac{p_xp_y}{2f}, ap_x^2+ bp_xp_y+ cp_y^2\right\}  \\ &=& \frac{1}{f}\left(p_x^3(fa_y) + p_x^2 p_y (fa_x + fb_y + 2 f_x a + f_y b)+ p_yp_x^2 (fb_x + fc_y+ f_x b + 2 f_y)+ p_y^3 (c_xf)\right),\end{array}\\$ i.e., is equivalent to the following system of PDE: 
 \begin{equation}\label{sys} 
 \left\{\begin{array}{rcc} a_y&=&0 \\ fa_x + fb_y + 2 f_x a + f_y b&=&0\\  
 fb_x + fc_y+ f_x b + 2 f_yc &=&0\\ c_x&=&0 \end{array}
 \right.\end{equation} 

 Thus, $a=a(x)$, $c=c(y)$, which is equivalent to 
  $B_1:= \frac{1}{\sqrt{|a|}} dx$ and $B_2:= \frac{1}{\sqrt{|c|}}dy$ are closed forms (assuming $a\ne 0$ and $c\ne 0$). \qed

\begin{rem} \label{rem3a} For further use let us formulate one more consequence of the equations \eqref{sys}: if $a\equiv c \equiv 0$ in a neighborhood of  a point, then $bf = \const$ implying     $F\equiv  \const \cdot H$ in the neighborhood.  
\end{rem}

Assume $a\ne 0$  ($c\ne 0$, respectively) 
 at a point $P_0$. For every point $P_1$ 
in a small neighbourhood $U$ of $P_0$ consider 
 \begin{equation} \label{normal}
 x_{new} :=\!\!\!\!\!\!\!\!\!
  \!\!\int\limits_{\begin{array}{c}
\gamma:[0,1]\to U \\ \gamma(0)=P_0 \\ \gamma(1)= P_1 \end{array}} 
\!\!\!\!\!\!\!\!\!\!\! \!\!\!\!\!\!\!B_1, \ \ \left(\textrm{$ y_{new} :=\!\!\!\!\!\!\!\!\!
  \!\!\int\limits_{\begin{array}{c}
\gamma:[0,1]\to U \\ \gamma(0)=P_0, \\  \gamma(1)= P_1 \end{array}} 
\!\!\!\!\!\!\!\!\!\!\! \!\!\!\!\!\!\!B_2, $  respectively}\right).\end{equation}

Locally, in the admissible coordinates, 
 the functions $x_{new}$ and $y_{new}$ are given by  
 \begin{equation*}
  x_{new}(x{})=\int_{x_0}^{x_{}} \frac{1}{\sqrt{|a(t)|} }\, dt,  \ \ \  \   y_{new}(y)=\int_{y_0}^{y_{}} \frac{1}{\sqrt{|c(t)|} }\, dt. 
\end{equation*}

The new  coordinates $(x_{new}, y_{new})$  \, $\bigl($or $(x_{new},y_{old})$ if $c_{old}\equiv 0$, 
or $(x_{old}, y_{new})$ if $a_{old}\equiv 0$$\bigr)$    are admissible.   In these coordinates,   the forms $B_1$ and $B_2$     are given by  $\textrm{sign}(a_{old}) dx_{new}$,  $\textrm{sign}(c_{old})dy_{new}$  (we assume $\textrm{sign}(0)=0$).

\subsubsection{ Proof of Theorem B}
We assume that $g$ of  signature $(+,-)$ on  $M^2$  
  admits  a nontrivial quadratic integral $F$ given by \eqref{integral1}. Consider 
 the matrix $
 F^{ij}= \begin{pmatrix} a & \tfrac{b}{2} \\ \tfrac{b}{2} & c \end{pmatrix} $. It 
  can be viewed as  a $(2,0)$-tensor: if we change the coordinate system and rewrite the  function $F$ in the new coordinates, the matrix changes according to the tensor rule. Then, 
  \begin{equation*}
   \sum_{\alpha=1}^2 g_{j\alpha} F^{i\alpha}
  \end{equation*}  
  is a $(1,1)$-tensor. 
   In  a neighborhood $U$  of almost every point the Jordan normal form of this $(1,1)$-tensor is one of the following matrices: 

 {Case 1} $\begin{pmatrix} \lambda & 0\\ 0 & \mu  \end{pmatrix} $, \ \  \  \  \  \  {Case 2}   $\begin{pmatrix} \lambda+ i \mu  & 0\\ 0 & \lambda- i \mu   \end{pmatrix} $, \ \   \  \  \  
  {Case 3} $\begin{pmatrix} \lambda  & 1\\ 0 & \lambda   \end{pmatrix} $, \\
  where $\lambda,  \mu:U\to \mathbb{R}$.  
 Moreover,  in view of Remark~\ref{rem3a}, there exists a neighborhood of almost every point such that $\lambda \ne \mu$ in Case 1 and $\mu \ne 0 $ in Case 2. 
  In the admissible coordinates, up to multiplication of $F $ by $-1$, and renaming $V_1\leftrightarrow V_2$,   
  Case 1 is equivalent to the condition 
    $a>0, \ c>0$,  Case 2 is equivalent to the condition   $a>0,$  $c<0$,  and Case 3 is equivalent to the   condition   $c\equiv 0$. 

  We now consider all three cases.

  \subsubsection{Case 1:  $a>0, \ c>0$. } 
 Consider the    coordinates \eqref{normal}. 
 In these coordinates,  $a=1$, $c=1$, and equations \eqref{sys} are: 

 \begin{equation*} 
 \left\{\begin{array}{rcc}  (fb)_y+ 2 f_x  &=&0,\\  
 (fb)_x  + 2 f_y&=&0. \end{array}
 \right.\end{equation*} 
 This system can be solved. Indeed, it is equivalent to 
 \begin{equation*} 
 \left\{\begin{array}{rcc}  (fb+ 2 f)_x  + (fb  + 2 f)_y&=&0,\\  
   (fb  -2 f)_x  -(fb- 2 f)_y&=&0, \end{array}
 \right.\end{equation*} 
 which,  after the change of cordinates $x_{new} = x+y$, $y_{new}= x-y$, has the form 
 \begin{equation*} 
 \left\{\begin{array}{rcc}  (fb+ 2 f)_x&=&0,\\  
   (fb  -2 f)_y&=&0, \end{array}
 \right.\end{equation*} 
 implying $fb+ 2f = Y(y)$, $fb-2f=X(x)$. Thus, 
 $f= \tfrac{Y(y)-X(x)}{4}$, $b=2 \tfrac{X(x)+Y(y)}{Y(y)-X(x)}$.

 Finally, in the new coordinates, the metric and the integral have (up to a  possible multiplication by a constant) the form

 \begin{equation*} 
 (X-Y)(dx^2 -dy^2)  
 \end{equation*}
 \and
 \begin{equation*}
 2\left(p_x^2 -  \tfrac{X(x)+Y(y)}{X(x)-Y(y)}(p_x^2-p_y^2) + p_y^2\right)= 4\frac{p_y^2 X(x) - p_x^2Y(y)}{X(x)-Y(y)}.  
 \end{equation*}
 Theorem B is proved under the assumptions of Case 1.

  \subsubsection{ Case 2:   
  $a>0$, $c<0$. }

 Consider the coordinates \eqref{normal}. 
 In these coordinates,  $a=1$, $c=-1$, and the equations \eqref{sys} are:

 \begin{equation*} 
 \left\{\begin{array}{rcc}  (fb)_y+ 2 f_x  &=&0,\\  
 (fb)_x  -2 f_y&=&0. \end{array}
 \right.\end{equation*} 
 We see that these conditions are the Cauchy-Riemann  conditions for the complex-valued 
 function $fb+ 2i \cdot f$. Thus, for an appropriate holomorphic 
 function $h= h(x+ i\cdot y)$,  we have $fb=\Re(h)$, $2f =\Im(h)$. 
 Finally, in a certain coordinate system the metric and the integral are 
 (up to multiplication by constants): 
  \begin{equation*}
 2\Im(h)dxdy \ \ \ \textrm{and} \ \ \ p_x^2 - p_y^2   + 2\tfrac{\Re(h)}{\Im(h)}p_xp_y.   
 \end{equation*}
 Theorem B is proved under the assumptions of Case 2.

 \subsubsection{Case 3: $a>0$, $c=0$.  } 

 Consider  admissible 
 coordinates $x, y$, such that  $x$ is 
 the coordinate from  \eqref{normal}. 
 In these coordinates,   $a=1$, $c=0$, and the equations \eqref{sys} are: 
 \begin{equation*} 
 \left\{\begin{array}{rcc}  (fb)_y+ 2 f_x  &=&0\\  
 (fb)_x  &=&0 \end{array}.
 \right.\end{equation*} 
 This system can be solved. Indeed, the second equation implies $fb= -Y(y)$. Substituting this in the first equation  we obtain 
 $Y'= 2f_x$ implying $$f= \frac{x}{2} Y'(y)+ \widehat{ Y}(y) \textrm{ \ \ and}  \ \ \ b= - \frac{Y(y)}{\frac{x}{2} Y'(y)+ \widehat{ Y}(y)}.$$ Finally, the metric and the integral are 
  \begin{equation}
  \label{answer:case3} 
 \left( \widehat{ Y}(y)+\frac{x}{2} Y'(y)\right)dxdy \ \ \ \textrm{and} \ \ \ p_x^2 - \frac{Y(y)}{\widehat{ Y}(y)+\frac{x}{2} Y'(y) }p_xp_y
 \end{equation}

 Moreover, by the change $y_{new}=\beta(y_{old})$ the metric and the integral  \eqref{answer:case3} will be  transformed to:
  \begin{equation*} 
 \left( \widehat{ Y}(y)\beta' +\frac{x}{2} Y'(y)\right)dxdy \ \ \ \textrm{and} \ \ \ p_x^2 +  \frac{Y(y)}{\widehat{ Y}(y)\beta'+\frac{x}{2} Y'(y) }p_xp_y
 \end{equation*}
 Thus, by putting $\beta(y) = \int_{y_0}^y \frac{1}{\widehat{ Y}(t)} dt$,  we can make the metric and the integral to be 
 $$  \left( 1+\frac{x}{2} Y'(y)\right)dxdy \ \ \ \textrm{and} \ \ \ p_x^2 -   \frac{Y(y)}{1+\frac{x}{2} Y'(y) }p_xp_y.
 $$
 Moreover, after the coordinate change $x_{new} =  \tfrac{x_{old}}{2}$ and dividing/multiplication of the metric/integral by ${2}$, the metric and the integral have  the form  from Theorem B
  \begin{equation*}
      \left( 1+ {x} Y'(y)\right)dxdy \ \ \ \textrm{and} \ \ \ p_x^2 -2  \frac{Y(y)}{1+ {x}{} Y'(y) }p_xp_y
 \end{equation*}

 Theorem B is proved.

  \end{appendix}

\end{document}